\documentclass[12pt]{amsart}
\usepackage{amssymb,eepic,epic}
\newtheorem{theorem}{Theorem}[section]
\newtheorem{corollary}[theorem]{Corollary}
\newtheorem{proposition}[theorem]{Proposition}

\newtheorem{lemma}[theorem]{Lemma}
\theoremstyle{definition}
\newtheorem{definition}[theorem]{Definition}
\theoremstyle{remark}
\newtheorem{remark}[theorem]{Remark}
\newtheorem{example}[theorem]{Example}
\newcommand\A{\mathcal{A}}
\newcommand\M{\mathcal{M}}
\newcommand\G{\mathcal{G}}

\renewcommand{\O}{\mathcal{O}}

\newcommand{\R}{\mathbb{R}}
\newcommand{\C}{\mathbb{C}}

\newcommand{\Z}{\mathbb{Z}}

\newcommand\lie[1]{\mathfrak{#1}}

\newcommand{\g}{\lie{g}}

\renewcommand{\t}{\lie{t}}
\newcommand{\Alc}{\lie{A}}
\renewcommand{\u}{\lie{u}}
\newcommand{\on}{\operatorname}

\newcommand{\length}{\on{length}}

\newcommand{\Ad}{ \on{Ad} } 
\newcommand{\Hol}{ \on{Hol} }

\newcommand{\Ind}{ \on{Ind}}
\renewcommand{\ker}{ \on{ker}}

\newcommand{\UU}{ \on{U}(1)}

\newcommand{\Vol}{  \on{Vol}}
\newcommand{\diag}{  \on{diag}}

\newcommand\dirac{/\kern-1.2ex\partial} 
\newcommand\qu{/\kern-.7ex/} 
\newcommand{\fus}{\circledast} \newcommand{\lev}{{\lambda}} 
\newcommand{\levi}{{k}} 
\newcommand{\LG}{\widehat{LG}} 
\newcommand{\Lg}{\widehat{L \g}} 
\newcommand{\Waff}{W_{\on{aff}}} 

\newcommand{\labell}\label

\newcommand{\lra}{\longrightarrow}
\newcommand{\hra}{\hookrightarrow}
\newcommand{\ra}{\rightarrow}

\renewcommand{\d}{{\mbox{d}}}
\newcommand{\ol}{\overline}
\newcommand\Gpar{\G_\partial}
\newcommand\Phinv{\Phi^{-1}}
\newcommand\phinv{\phi^{-1}}

\newcommand\lam{\lambda}

\newcommand\Sig{\Sigma}
\newcommand\sig{\sigma}
\newcommand\eps{\epsilon}
\newcommand\Om{\Omega}
\newcommand\om{\omega}

\newcommand{\del}{\delta}
\newcommand{\f}{\frac}
\renewcommand{\c}{\mathcal{C}}

\newcommand{\p}{\partial}
\renewcommand{\l}{\langle}
\renewcommand{\r}{\rangle}

\newcommand{\maxx}{0}

\begin{document}

\title[Fusion and cobordism]
{Fusion of Hamiltonian loop group manifolds\\ and cobordism}

%
\author{E. Meinrenken} \address{Massachusetts Institute of Technology,
Department of Mathematics, Cambridge, Massachusetts 02139}
\thanks{Supported by a Feodor Lynen fellowship from the Humboldt foundation.}
\email{mein@math.mit.edu}
\author{C. Woodward} \address{Harvard University, Department of
Mathematics, 1 Oxford Street, Cambridge, Massachusetts 02138}
\thanks{Supported by an NSF Postdoctoral Fellowship.}  
\date{July 24, 1997}

\email{ woodward@math.harvard.edu}

\begin{abstract}  
We construct an oriented cobordism between moduli spaces of flat
connections on the three holed sphere and disjoint unions of toric
varieties, together with a closed two-form which restricts to the
symplectic forms on the ends.  As applications, we obtain formulas for
mixed Pontrjagin numbers and Witten's formulas for symplectic volumes.
\end{abstract}

\maketitle
\tableofcontents


\section{Introduction}

Let $G$ be a compact, connected, simply connected, simple 
Lie group and $\Sig$ a compact
oriented 2-manifold, possibly with boundary.  A standard approach to
computing invariants of the moduli spaces of flat $G$-connections over
$\Sig$ is to study the behavior of these invariants when two boundary
circles are glued together.  Since any $\Sig$ of genus at least two
admits a pants decomposition, this approach reduces the computation of
these invariants to that of the moduli spaces associated to pairs of
pants (three-holed spheres).

In this paper, we develop a cobordism approach to computing invariants
of these moduli spaces.  We apply our technique to compute  
mixed Pontrjagin numbers and symplectic volumes.
In a sequel \cite{me:co2} to this paper, 
we apply our method to compute the
coefficients of the fusion ring (Verlinde algebra).

The main result is as follows.  
Let $T\subset G$ be a maximal torus,  
$\t_+\subset \t$ a choice
of a closed positive Weyl chamber and $\Alc\subset\t_+$ the closed
fundamental Weyl alcove.
We interpret $\Alc$ as the set of
conjugacy classes in $G$ since for every conjugacy class $\c \subset
G$ there is a unique $\mu\in\Alc$ with $\exp(\mu)\in \c$.
For $\mu_1,\mu_2,\mu_3\in\Alc$ let 
$$ \M(\Sig_0^3,\mu_1,\mu_2,\mu_3)$$
denote the moduli space of flat $G$-connections on the three-holed
sphere $\Sig_0^3$ for which the holonomies around the three boundary
components lie in the conjugacy classes labeled by the $\mu_j$.
A result of L. Jeffrey \cite{je:em} identifies the moduli
spaces $\M(\Sig_0^3,\mu_1,\mu_2,\mu_3)$ when the holonomies are small.
Use the normalized invariant inner product on $\g$ to identify
$\g\cong\g^*$ and let $\O_{\mu_j}$ be the coadjoint orbit through
$\mu_j$, equipped with the Kirillov-Kostant-Souriau symplectic form.
Let $*:\ \t_+ \ra \t_+$ denote the involution defined by
$*\mu=-w_0\cdot\mu$ where $w_0$ is the longest Weyl group element.
Jeffrey proves in \cite{je:em} that for $\mu_1,\mu_2,\mu_3\in\Alc$
sufficiently small, there is a symplectomorphism
$$ \M(\Sig_0^3,\mu_1,\mu_2,\mu_3)\cong (\O_{* \mu_1} \times \O_{*
\mu_2} \times \O_{* \mu_3}) \qu G $$
of the moduli space with a symplectic reduction of a triple product of
coadjoint orbits under the triagonal action.
The main theorem of this paper is the following result for more general
holonomies:
\begin{theorem} \label{MainTheorem}\label{MainThm}
For $\mu_1,\mu_2,\mu_3\in \on{int}(\Alc)$ generic  there is an
oriented orbifold cobordism
$$
\M(\Sig_0^3,\mu_1,\mu_2,\mu_3)\sim \coprod_{w\in\Waff^+}
(-1)^{\length(w)}\, (\O_{* \mu_1} \times \O_{* \mu_2} \times \O_{w\,* \mu_3})
\qu G.
$$
Here the signs indicate a change in orientation relative to the
symplectic orientation, and $\Waff^+$ is the set of all $w$ in the
affine Weyl group $\Waff$ such that $w\Alc\subset\t_+$.  The
symplectic forms extend to a closed 2-form over the cobordism. For
$G=\on{SU}(n)$ both sides are smooth manifolds and the cobordism is a
manifold cobordism.
\end{theorem}  
Cobordisms of orbifolds with closed 2-forms were introduced and 
studied by Ginzburg-Guillemin-Karshon in \cite{gi:co}.
The genericity assumption in Theorem \ref{MainThm}
guarantees that both sides have at worst
orbifold singularities.  In fact, the cobordism can be carried further
so that the moduli space is cobordant to a disjoint union of toric
varieties (Theorem \ref{ToricVarieties}).  An application is the
computation of the symplectic volume and mixed Pontrjagin numbers of
the moduli space of the three-holed sphere (Theorem \ref{VolPontr}).
Just as for toric varieties, the mixed Pontrjagin numbers are obtained
by applying differential operators to the volume functions. 
In \cite{me:wi} we obtain
Witten's volume formulas for arbitrary compact, oriented 2-manifolds
(cf. Theorem \ref{WitThm}) by gluing.  These formulas were proved in most
cases by Witten \cite{wi:qg}.  Alternative proofs and extensions were
given by Liu \cite{li:hk,li:h2}, Jeffrey-Kirwan \cite{je:in} and
Jeffrey-Weitsman \cite{jw:va}.

The idea of our proof is based on the ``classical analog'' of the
principle (proved e.g. in \cite{te:lc}) that ``induction from
$G$-representations to representations of the loop group $LG$ takes
tensor product to fusion product''.  We introduce a notion of fusion
product for Hamiltonian loop group manifolds, which preserves
properness and level of the moment map.  The moduli space
$\M(\Sig^3_0, \mu_1,\mu_2,\cdot)$ of flat connections on the
three-holed sphere with fixed holonomy around two boundary components
is the fusion product of two coadjoint orbits of $LG$.  We also
define a notion of induction $\Ind$ of Hamiltonian $G$-manifolds to
Hamiltonian $LG$-manifolds.  We show that the induction of a direct
product is cobordant (by a topologically trivial cobordism) to the
fusion product of the inductions.  In particular, we find an 
equivariant topologically trivial cobordism of Hamiltonian $LG$-manifolds,
$$ \M(\Sig^3_0,
\mu_1,\mu_2,\cdot)\sim \Ind(\O_{*\mu_1}\times \O_{*\mu_2}).$$
An application of the ``cobordism
commutes with reduction'' principle of \cite{gi:co} 
yields Theorem \ref{MainThm}.

We give a brief outline of the contents.  In Section 2 we review the
definition of a Hamiltonian group action, and the definition of
Hamiltonian cobordism according to Ginzburg-Guillemin-Karshon
\cite{gi:co}.  The orbifolds appearing in Theorem \ref{MainTheorem}
are naturally viewed as quotients of Hamiltonian $LG$-manifolds, which
are discussed in Section 3.  In Section 4 we introduce the notion of
fusion product of Hamiltonian $LG$-manifolds, and prove Theorem
\ref{MainThm}.  In Section 5 we present the application of our results
to mixed Pontrjagin numbers and symplectic volumes of moduli spaces.
A key technical ingredient in our proof of Theorem \ref{MainTheorem}
is a Duistermaat-Heckman principle for Hamiltonian loop group
manifolds which we prove in Section 6.

\section{Preliminaries}

In this section we will review the material we want to generalize to
the $LG$-equivariant setting.

\subsection{Hamiltonian G-manifolds}

Let $G$ be a compact Lie group acting on a manifold $M$.  We denote by
$T$ a choice of maximal torus, by $\t^*_+\subset\t^*\subset\g^*$ a
choice of positive Weyl chamber and by $\Lambda^*_+\subset\t^*_+$ the
corresponding set of dominant weights.  For any $\xi$ in the Lie
algebra $\g$ of $G$ we have a generating vector field $\xi_M \in
\on{Vect}(M)$ whose flow is given by $ m \mapsto \exp(t\xi)m$.  Let
$\om \in \Omega^2(M)$ be a closed $G$-invariant $2$-form.  The action
of $G$ on $(M,\om)$ is called Hamiltonian if there exists a {\em
moment map} $\Phi:\,M\to \g^*$ which is equivariant with respect to
the coadjoint action on $\g^*$, and satisfies
\begin{equation} \iota(\xi_M)\omega = \d\l\Phi,\xi\r. 
\labell{DefnMomentMap}\end{equation}
Note that we do not require the two-form $\om$ to be be
non-degenerate (i.e. symplectic).  

An important consequence of this definition is that for any 
$m\in M$, the stabilizer algebra $\g_m$ annihilates the image 
of the tangent map $\d_m\Phi$:
$$ 
\g_m\subset \on{im}(\d_m\Phi)^0
$$
(with equality if $M$ is symplectic). 
In particular, if $\mu$ is a regular value of the moment map then the
action of $G_\mu$ is locally free at points of the submanifold
$\Phi^{-1}(\mu)$.  The quotient (reduced space at $\mu$)
$$M_\mu=\Phi^{-1}(\mu)/G_\mu$$ 
is therefore an orbifold, 
and carries a canonical closed 2-form $\omega_\mu$. For $\mu=0$ we 
will also use the notation 
$$ M_0=M\qu G=\Phinv(0)/G, $$  
this is convenient if there is another group action on $M$.
We say that $(M,\om,\Phi)$ is {\em pre-quantizable}
if there exists a $G$-equivariant Hermitian line bundle $L$
with invariant connection $\nabla$ such that 
\begin{equation}\labell{Form1}\omega=
{\textstyle \f{i}{2\pi}}\on{curv}(\nabla),\end{equation} 
\begin{equation}\labell{Form2}
2\pi i\,\l\Phi,\,\xi\r=\on{Vert}(\xi_L)\in C^\infty(M,\on{End}(L)).
\end{equation}
Here we have identified sections of the bundle $\on{End}(L)$ of
Hermitian bundle endomorphisms with imaginary-valued functions
$C^\infty(M,i\R)$;  $\xi_L$ denotes the fundamental vector
on the total space of $L$ and $\on{Vert}:\,TL\to TL$ the vertical
projection given by the connection. $(L,\nabla)$
is called a {\em pre-quantum line bundle} for $(M,\om,\Phi)$.  
Clearly a necessary condition for the existence of $(L,\nabla)$
is that $\om \in H^2(M,\Z)$.  If $G$ is connected, simply connected
and compact, this condition is also sufficient.  
Conversely, any $G$-equivariant Hermitian line bundle with invariant
connection gives $M$ the structure of a Hamiltonian $G$-manifold
via the equations (\ref{Form1}),(\ref{Form2}).

If $0$ is a regular value of $\Phi$, the reduced space $M\qu G$ has 
pre-quantum line bundle 
$$ L\qu G= L_0=(L|\Phinv(0))/G $$ 
(in the case that the $G$-action on $\Phinv(0)$ is not free, this is an
orbi-bundle). More generally, if $\mu\in\Lambda^*_+$ is a weight for 
$G$, let $\C_\mu$ denote the 1-dimensional $G_\mu$-representation 
defined by $\mu$. If $\mu$ is a regular value for $\Phi$, 
a pre-quantum (orbi)-bundle for $M_\mu$ is given by 
$$ L_\mu^{shift}=\big((L\otimes \C_{*\mu})|\Phinv(\mu)\big)/ G_\mu. $$

\subsection{Cobordism}
We recall the notion of {\em cobordism} in the category of Hamiltonian
$G$-manifolds, as introduced by Guillemin-Ginzburg-Karshon \cite{gi:co}.
Let $G$ be a compact Lie group.

\begin{definition}[Hamiltonian Cobordism]
\labell{DefinitionCobordism} \cite{gi:co}
A cobordism between two oriented Hamiltonian $G$-manifolds
$(M_1,\om_1,\Phi_1)$ and $(M_2,\om_2,\Phi_2)$ with proper moment maps
is an oriented Hamiltonian $G$-manifold with boundary $(N,\om,\Phi)$
with proper moment map such that $\p N=M_1\cup (-M_2)$ and such that
$\om$ resp. $\Phi$ pull back to $\om_i$ resp. $\Phi_i$. We write
$M_1\sim M_2$, or sometimes $(M_1,\om_1,\Phi_1)\sim
(M_2,\om_2,\Phi_2)$.
\end{definition}

This notion of cobordism has the drawback that it is not well-behaved
under reduction, since reduced spaces of a Hamiltonian $G$-manifold
are generically orbifolds. We will use the terminology ``orbifold
cobordism'' if $M_1,\,M_2,\,N$ in the above definition are allowed to 
have orbifold singularities. 
(See Druschel \cite{dr:oo} for information on oriented orbifold cobordisms.)

The following Lemma is completely obvious, but has a lot of interesting 
consequences as shown in \cite{gi:co}. 

\begin{lemma}
\labell{CobComRed}
Suppose $(N,\om,\Phi)$ gives a cobordism 
$(M_1,\om_1,\Phi_1)\sim (M_2,\om_2,\Phi_2)$. 
If $\mu$ is a regular value of $\Phi$, the reduced space 
$N_\mu$ gives an orbifold cobordism
$((M_1)_\mu,(\om_1)_\mu)\sim \big((M_2)_\mu,(\om_2)_\mu\big)$.  
\end{lemma}

\begin{example} \labell{Product}
Suppose that $M_1=M_2=M$ are compact oriented $G$-manifolds and that the
equivariant closed 2-forms $\om_i+2\pi i\Phi_i$ are cohomologous,
i.e. $\exists\beta\in\Om^1(M)^G$ such that
$$\om_2-\om_1=\d\beta,\,\, \Phi_2-\Phi_1=-\beta^\sharp$$
where $\beta^\sharp:\,M\to \g^*$ is defined by $\l\beta^\sharp,\xi\r
=\iota(\xi_M)\beta$. 
Let 
$N=M\times [0,1]$ with coordinates $(m,t)$, and let
$$\om=\om_1+\d(t\beta),\,\,\Phi(m,t)=(1-t)\Phi_1(m)+t\Phi_2(m).$$ 
Then $(N,\om,\Phi)$ provides a cobordism of Hamiltonian $G$-manifolds,
$$(M_1,\om_1,\Phi_1)\sim (M_2,\om_2,\Phi_2).$$
From this ``trivial''
cobordism, non-trivial examples are obtained by reduction. We remark
that it is not sufficient in this example to assume properness of the
moment maps $\Phi_i$ since this does not imply properness of $\Phi$ in
general.
\end{example}

\section{Hamiltonian actions of loop groups}

\subsection{Loop groups}
Let $G$ be a compact, connected, simply connected Lie group. Let the 
Lie algebra $\g$ be equipped with the invariant inner product, 
normalized  by the requirement that on each simple factor, the long 
roots have length $\sqrt{2}$. 

We define the loop
group $LG$ as the Banach Lie group of maps $S^1 \ra G$ of some fixed
Sobolev class $s > 1/2$. $LG$ is a semi-direct product $LG=\Omega
G\rtimes G$ where the subgroup $\Omega G$ of based loops is the kernel
of the evaluation map $LG\to G,\,g\mapsto g(1)$ and $G$ acts on $\Om
G$ by conjugation.  Let
$$ 1 \ra \UU \ra \LG \ra LG \ra 1 $$
denote the basic central extension of $LG$ (see \cite{ps:lg}). The 
corresponding Lie algebra extension $\widehat{L\g}$  
is given by $L\g\times \R$, with bracket   
$$ [(\xi_1,t_1),(\xi_2,t_2)]=
\big([\xi_1,\xi_2],\,\oint \xi_1\cdot\,
\d \xi_2\big). $$ 
using the normalized inner product on $\g$. The inner product on 
$\g$ gives defines a natural $LG$-invariant $L^2$-metric on $L\g$. 

\subsubsection{The affine coadjoint action}
Let $L\g^*$ denote the space of $\g$-valued $1$-forms
$\Omega^1(S^1,\g)$ of Sobolev class $s-1$.  We consider $L\g^*$ as a
subset of the topological dual space to $L\g$ via the natural pairing
given by integration.  The coadjoint action of $LG$ on
$\widehat{L\g}^*=L\g^*\times\R$ is given by
$$ g\cdot (\mu,\lambda)=(\Ad_g(\mu)-\lambda\, \d g\,g^{-1},\,\lambda).$$
By the {\em affine action at level $\lev$} we mean the $LG$-action on
$L\g^*$ corresponding to the identification with the hyperplane
$L\g^*\times \{\lev\}\subset \widehat{L\g}^*$. Clearly, the actions at
non-zero level are equivalent up to rescaling, so that we will always
consider the level $\lev=1$ unless specified otherwise.  
We will consider $\g$ as a subset of $L \g^*$ by the embedding
$\xi \mapsto \xi \d \theta / 2 \pi $.

\subsubsection{The holonomy map}

The action of $LG$ on $L\g^*$ at level $1$ can be identified with the
action of the gauge group of $S^1$ on connections on the trivial
principal $G$-bundle over $S^1$.  Taking the holonomy of such a
connection around $S^1$ we obtain a smooth submersion
$$ \Hol:\,L\g^*\to G $$
with the equivariance property $ \Hol(g\cdot\mu)=\Ad_{g(1)}\Hol(\mu)$.
The restriction of $\Hol$ to $\g \subset L \g^*$ is the exponential
mapping.  For any point $\mu\in L\g^*$ the evaluation mapping $LG\to
G,\, g\mapsto g(1)$ gives an isomorphism
\begin{equation}  (LG)_\mu \cong Z_{\Hol(\mu)}
\label{HolonomyIsomorphism}
\end{equation} 
where $ Z_{\Hol(\mu)}$ is the centralizer of the holonomy of $\mu$
(see \cite{ps:lg}). In particular, the coadjoint action of the subgroup
$\Om G$ of based loops is free, and the holonomy map may be viewed as
the quotient map $L\g^*\mapsto L\g^*/\Om G=G$.

\subsubsection{The inversion map}
There is a natural automorphism, the {\em inversion map} 
\begin{equation}\labell{Inversion}
I^*:\,\widehat{LG}\to\widehat{LG}
\end{equation}
which is defined on the Lie algebra level as follows.
Let $I:\,S^1\to S^1$ denote the map $z\mapsto z^{-1}$. 
The pullback action on $L\g$ extends to a Lie algebra homomorphism 
$I^*:\,\widehat{L\g}\to\widehat{L\g}$ by  $I^*(\xi,t)=(I^*\xi,-t). $
Indeed, since $I$ reverses the orientation on $S^1$ the cocycle 
$c(\xi_1,\xi_2)=\oint \xi_1\d\xi_2$ transforms according 
to $c(I^*\xi_1,I^*\xi_2)=-c(\xi_1,\xi_2)$.
Note that dual action on $\widehat{L\g}^*$ is given by 
$$I^*(\mu,\,\lev)=(-I^*\mu,-\lev)$$
which in particular changes the sign of levels.

\subsubsection{The involution $*$}

Recall that the involution $*:\,\t_+\to \t_+$ extends to a Lie algebra
homomorphism $*:\,\g\to \g$ which in turn exponentiates to a map
$*:\,G\to G$. In a faithful matrix representation for $G$, this is the
map $g\mapsto (g^{-1})^t$ where $t$ denotes ``transpose''. Composition
with $*:\,G\to G$ transforms any unitary $G$-representation into the
dual representation.  The induced homomorphism $*:\,LG\to LG$ extends
to a homomorphism of $\widehat{LG}$, given on the Lie algebra by
$$ *:\,\widehat{L\g}\to \widehat{L\g},\,(\xi,t)\mapsto (*\xi,t). $$
Notice that the dual map $*$ on $\widehat{L\g}^*$ is level preserving. 
Note also that if $\mu\in L\g^*$ is at level 1, $\Hol(*\mu)=*\Hol(\mu)$
since $\Hol$ is equivariant and restricts to the exponential map
on $\g\subset L\g^*$.

\subsection{Hamiltonian LG-manifolds}

Let $M$ be a Banach manifold, together with a closed 2-form $\omega$.
We call $\omega$ symplectic if $\omega$ is weakly non-degenerate, that
is, if $\omega$ defines an injection $T_m M \ra T_m^* M$.  An action
of $\LG$ on $M$ which preserves $\omega$ is called {\em Hamiltonian}
if there exists a moment map $\Phi:\,M\to \Lg^*$ satisfying
\eqref{DefnMomentMap}. Since the action of $LG$ on $\widehat{L\g}^*$
has fixed point set $\{0\}$, the equivariance condition implies that
$\Phi$ is unique if it exists.  An $\LG$-equivariant Hermitian line
bundle $L\to M$ with connection $\nabla$ satisfying \eqref{Form1} and
\eqref{Form2} is called a pre-quantum line bundle.

We call $M$ a Hamiltonian $LG$-manifold at level $\lev$ if the central
$\UU$ acts trivially with moment map $\lev$. If $M$ admits a pre-quantum line 
bundle $L$ this implies $\lev=k\in\Z$, where $k$ is the weight for the 
action of the central $\UU$ on $L$.

Note that if $(M,\omega,\Phi)$ is a Hamiltonian $LG$-manifold at non-zero
level $\lev$ then $(M,\f{\omega}{\lev},\f{\Phi}{\lev})$ is a Hamiltonian
$LG$-manifold at level $1$. We will therefore always assume $\lev=1$
unless specified otherwise.

\subsubsection{The holonomy manifold}
Notice that by the equivariance condition, the action of the subgroup
of based loops $\Om G $ is {\em free}.  Let $\Hol(M):=M/\Om G$ and
$\Hol(\Phi):\Hol(M)\to G$ the map induced by $\Phi$.  We have a
commutative diagram
$$ \begin{array}{ccc}
M & \stackrel{\Phi}{\lra} & L\g^* \\
\downarrow & & \downarrow \\
\Hol(M) & \stackrel{\Hol(\Phi)}{\lra} & G 
\end{array} $$
There is a unique manifold structure on the {\em holonomy manifold}
$\Hol(M)$ such that the quotient map $M\to \Hol(M)$ is a
submersion. (Existence follows e.g. from the existence of slices, or
from the discussion in Section \ref{CrossSections} below.  For
uniqueness see \cite{ab:ma}, Prop. 3.5.20.) The map $\Hol(\Phi)$ is
smooth and equivariant.
  
The $LG$-manifold $M$ together with the moment map $\Phi$ may thus be
viewed as the pull-back of the universal $\Om G$-principal bundle
$E\Om G= L\g^*\to B \Om G=G$ under the map $\Hol(\Phi)$.  Note that
$\Hol(M)$ is compact (in particular finite-dimensional) if and only if
$\Phi$ is proper. The 2-form on $M$ is not basic and therefore does not
descend to $\Hol(M)$. However, as shown in \cite{al:mom}
there exists a canonical 2-form 
$\varpi$ on $L\g^*$ such that $\om+\Phi^*\varpi$ is basic and 
descends to 
a 2-form on $\Hol(M)$ for which $\Hol(\Phi)$ may be interpreted as 
a group valued moment map. See \cite{al:mom} for a detailed discussion. 

\subsubsection{The involution $M\mapsto M^*$} 

Suppose $(M,\omega,\Phi)$ is a Hamiltonian $LG$-manifold at level $1$.
Define a new action on $(M,\omega)$ by composing the action map $LG\to
\on{Diff}(M)$ with the involution $*:\,LG\to LG$.  The new action is
once again Hamiltonian at level $1$, and the moment map is
$\Phi^*:=*\Phi$. We denote the resulting space by $M^*$.  Clearly
$\Hol(M^*)$ is just $\Hol(M)$, with new $G$-action defined by
composing with $*:\,G\to G$, and moment map $\Hol(\Phi)^*=*\Hol(\Phi)$.

\subsubsection{Line bundles}
An $\LG$-equivariant line bundle $L$ over an $LG$-Banach manifold $M$
is called at level $k\in\Z$ if the central circle acts with weight
$k$.  If $k=0$, this means that $L$ is an $LG$-equivariant line bundle
which descends to a $G$-equivariant line bundle $\Hol(L) := L/ \Omega
G$ over $\Hol(M)$.  It follows that the isomorphism classes of (level
$0$) $LG$-equivariant line bundles $L\to M$ are classified by the
equivariant cohomology $H^2_G(\Hol(M),\Z)$.

\subsection{Coadjoint orbits}

Coadjoint $LG$-orbits provide basic examples for Hamiltonian
$LG$-manifolds, the moment map being simply the embedding. In
particular, the coadjoint orbit $LG\cdot 0$ is $\Omega G=LG/G$.  The
isomorphism (\ref{HolonomyIsomorphism}) shows that every stabilizer
group $(LG)_\mu$ is compact and connected (because for a compact
simply connected Lie group, the centralizer of any element is
connected).  Let $G$ be a simple, simply-connected, compact,
connected Lie group, with fundamental alcove $\Alc$.  Every coadjoint
$LG$-orbit passes through exactly one point of $\Alc$, considered as a
subset of $L\g^*$ via the embeddings
$$\Alc\to \t \to \g \to L\g^* .$$ 
Using the exponential mapping we can view $\Alc$ also as a subset of
$T\subset G$, and every $\Ad(G)$-orbit passes through exactly one
point of $\Alc$.  We therefore obtain a series of identifications
$$ \Alc\cong T/W\cong G/\Ad(G)\cong \t/\Waff\cong L\g^*/LG. $$ 
One also finds that for any open face $\sig\subset\Alc$, the
stabilizer groups $(LG)_\mu$ resp. centralizers $Z_{\exp\mu}$ for all
points $\mu\in\sig$ are the same, they will therefore be denoted by
$(LG)_\sig$ resp. $Z_\sig$.  The Dynkin diagrams of the semi-simple
parts of these groups are obtained from the extended Dynkin diagram of
$G$ by deleting the vertices corresponding to affine simple roots not
vanishing on $\sigma$. In particular, the groups $(LG)_\sig$
corresponding to the vertices of $\Alc$ are semi-simple (but not
always simply-connected).  There is a natural partial ordering on the
open faces of $\Alc$ defined by inclusion: We write $\sig \prec \tau$
if $\sig$ is properly contained in $\ol{\tau}$. Let
$\alpha_1,\ldots,\alpha_l$ be the simple roots of $G$, and
$\alpha_{\maxx}$ the highest root.  We fix an orientation on $G$, and
let $T$ have the orientation induced by the choice of positive roots.

\begin{lemma}(Properties of the stabilizer groups $(LG)_\sig$).
\labell{PropertiesStabilizers}
\begin{enumerate}
\item
Every $(LG)_\sig$ contains the maximal torus $T$. Define the
fundamental Weyl chamber for $(LG)_\sig$ as the cone over $(\Alc-\mu)$
where $\mu\in\sig$.  Together with the orientation on $T$ this induces
an orientation on $(LG)_\sig$. The simple roots for $(LG)_\sig$ are
precisely those roots $\alpha$ in the collection
$\{\alpha_1,\ldots,\alpha_l,-\alpha_{\maxx}\}$ such that $\alpha$
vanishes on the span of $\sig-\mu$.
\item
For $ \sig\prec \tau$, one has $ (LG)_\sig\supset (LG)_\tau$. The
quotient $(LG)_\sig/(LG)_\tau$ carries a canonical invariant complex
structure.  
\item 
 If $0\in\ol{\sig}$ then $(LG)_\sig$ is contained in the  subgroup 
$G\subset LG$ of constant loops.
\item
For every $\sig$, the Lie algebra $(L\g)_\sig$ of $(LG)_\sig$ 
has a unique $(LG)_\sig$-invariant complement in $L\g$, given 
by the $L^2$-orthogonal complement $(L\g)_\sig^\perp$.  
\end{enumerate}
\end{lemma}
If $G$ is a simply-connected, compact, connected Lie group, and $G =
G_1 \times \ldots \times G_r$ its decomposition into simple factors,
then coadjoint orbits of $LG = LG_1 \times \ldots \times LG_r$ are
products of coadjoint orbits of the $LG_i$.  We denote by $\Alc$ the
product of the fundamental alcoves $\Alc_i$ for the $G_i$.

To describe the Kirillov-Kostant-Souriau
(KKS) form on coadjoint orbits $(LG)\cdot\mu$ 
let $\del_\mu$ denote the elliptic operator
$$ \del_\mu: L\g\to L\g^*,\, \zeta\to \d \zeta+\,[\mu,\zeta]. $$
Then $\del_\mu(\xi)$ is the value at $\mu$ for the fundamental vector
field defined by $\xi$. In particular, the kernel of $\del_\mu$ is
equal to $(L\g)_\mu$ and the image is equal to the tangent space to
the orbit $T_\mu(LG\cdot\mu)$. The KKS-form $\nu_\mu\in\Om^2(LG\cdot\mu)$ 
is the unique invariant 2-form given on the tangent space 
$T_\mu(LG\cdot\mu)$ by 
$$
\nu_\mu(\del_\mu(\xi_1),\del_\mu(\xi_2))=
\l (\mu,1),[(\xi_1,0),(\xi_2,0)] \r =
\l \mu,[\xi_1,\xi_2]\r +\oint
\xi_1\cdot \d\xi_2
$$
This can be rewritten
$$
\nu_\mu(\del_\mu(\xi_1),\del_\mu(\xi_2))=
\oint
\xi_1\,\cdot\, \del_\mu (\xi_2).
$$

\begin{remark}
The coadjoint orbits $LG\cdot\mu$ for $\mu\in\Alc$ 
have a canonical invariant complex structure, making them into 
$LG$-K\"ahler manifolds. Using the operator $\del_\mu$ the 
complex structure can be described as follows. Identify the 
tangent space $T_\mu(LG\cdot\mu)$ with the $L^2$-orthogonal 
complement $\ker(\del_\mu)^\perp$. 
Consider $*\del_\mu$ (where $*$ is the Hodge operator) as a
first order skew adjoint 
differential operator on $C^\infty(S^1,\g)$. 
The square $(*\del_\mu)^2$ is a
non-positive operator, so that $|*\del_\mu|:=\sqrt{-(*\del_\mu)^2}$ is a
well-defined first order pseudo-differential operator and its Greens
operator $|*\del_\mu|^{-1}$, defined to be the identity on the kernel
of $\del_\mu$, has order $-1$. The complex structure is given by 
the zeroth order pseudo-differential operator 
$$
J_\mu:=|*\del_\mu|^{-1}\,\, (*\del_\mu)
$$
on $\ker(\del_\mu)^\perp\subset L\g$. 
The proof for the fundamental homogeneous space $LG\cdot 0=LG/G$ 
can be found in Freed \cite{fr:lg}; 
the extension to the general case is straightforward.
\end{remark}

\begin{remark}
For any coadjoint orbit $LG\cdot\mu$, we have $(LG\cdot\mu)^*=
LG\cdot *\mu$. 
\end{remark}

\begin{remark}
The holonomy manifold of a coadjoint orbit $LG\cdot\mu$ is just 
the conjugacy class $G\cdot \exp(\mu)$, and the map 
$\Hol(\Phi)$ is the embedding into $G$.
\end{remark}

\subsection{Moduli spaces of flat connections} 
\labell{Fusion}

Let $\Sig$ be a compact oriented 2-manifold with $b$ boundary
components.  We denote by $\iota:\,\partial\Sigma\hra \Sigma$ the
inclusion of the boundary.  Fix $s>1$, and let
$\A(\Sig)\cong\Omega^1(\Sig,\g)$ the space of connections of Sobolev
class $s-\f{1}{2}$ in the trivial principal bundle $\Sig\times G$, and
$\G(\Sig)$ the gauge group, consisting of maps $\Sig\to G$ of Sobolev
class $s+\f{1}{2}$.  
Let $\iota:\,\p\Sig\hra \Sig$ be the inclusion of the boundary. 
Recall that restriction $\iota^*$ to the boundary results
in the loss of half a derivative, so that there is a continuous map
$\G(\Sig)\to \G(\p\Sig)\cong LG^b$ which is surjective since $\pi_1(G)
= 0$.  The kernel $\G_\p(\Sig)$ consists of gauge transformations that
are the identity on the boundary.  
According to Atiyah-Bott \cite{at:mo},
the gauge group action on the space $\A(\Sig)$ with symplectic form
$$ \omega_A(a_1,a_2)= 
\int_\Sigma a_1 \stackrel{ _\cdot}{\wedge} a_2
\,\,\,\,\,\,\,\,\,(a_i\in T_A\A(\Sig)\cong\Omega^1(\Sig,\g)) $$
is Hamiltonian, with moment map given by 
$$ 
\Psi:\,\A(\Sig)\to\Om^2(\Sig,\g)\oplus \Om^1(\p\Sig,\g),\ 
A\mapsto (\on{curv}(A),\,\iota^* A), 
$$
that is, 
$$ 
\l\Psi(A),\xi\r= \int_\Sigma \on{curv}(A) \cdot\xi + \int_{\partial
\Sigma} \iota^*\,(A\cdot\xi).
$$
Here $\on{curv}(A)$ denotes the curvature of $A$, and 
we have chosen the orientation on the boundary $\p \Sig$ to be 
minus the
orientation induced from $\Sig$.  The moment map for the action of
$\Gpar(\Sig)$ on $\A(\Sig)$ can be identified with $A \mapsto
\on{curv}(A)$ and hence the symplectic quotient of $\A(\Sig)$ by
$\Gpar(\Sig)$ is
$$ \M(\Sig) := \A_F(\Sig)/\Gpar(\Sig) $$
where $\A_F(\Sig) \subset \A(\Sig)$ is the space of flat connections.
If $\partial \Sig=\emptyset$ then $\M(\Sig)$ is a compact, finite
dimensional stratified symplectic space (in general singular).  On the
other hand, if $\partial \Sig\not=\emptyset$ then according to
Donaldson \cite{do:bv} $\M(\Sig)$ is a smooth infinite-dimensional
symplectic manifold.  It has a residual Hamiltonian action of
$\G(\partial \Sig)\cong LG^b$ with moment map
$$ \Phi: \M(\Sig) \ra \Omega^1(\partial \Sig,\g), \ \ [ A ] \mapsto 
\iota^* A$$
and carries a natural $\widehat{LG^b}$-equivariant pre-quantum line bundle 
$L(\Sig)$ with connection. 
\begin{remark}\label{Selfadjoint}
For any 2-manifold $\Sig$, the symplectomorphism 
$\A(\Sig)\to \A(\Sig),\,A\mapsto A^*$ induced by 
the involution $*:\,\g\to \g$ preserves $\A_F(\Sig)$ and 
therefore
descends to a symplectomorphism 
$\M(\Sig)^*\cong\M(\Sig)$.
\end{remark}
\begin{remark}
The holonomy manifold $\Hol(\M(\Sig))$ may be interpreted as the 
quotient of the space $\A_F(\Sig)$ of flat connections by 
those gauge transformations which are the identity on given 
base points of the boundary circles.
\end{remark}

\begin{remark}
For the sake of completeness we recall from \cite{se:ln,do:bv} the
holomorphic description of $\M(\Sig)$.  The moduli space $\M(\Sig)$
inherits a complex structure from a choice of complex structure on
$\Sig$.  The action of $\G(\p\Sig)=LG^b$ complexifies to an action of
the $b$-fold complex loop group $\G_\C(\partial \Sig)=LG_\C^b$ and
(for $b> 0$) $\M(\Sig)$ is a homogeneous space
$$ \M(\Sig) = \G_\C(\partial \Sig)/\G_\C(\Sig) $$
where $\G_\C(\Sig)$ is the group of holomorphic maps of $\Sig$
into $G_\C$.  
\end{remark}

\begin{example}\labell{ex:2Holes}
\begin{enumerate}
\item \labell{1Hole} The moduli space $\M(\Sig^1_0)$ for the disk is
the fundamental homogeneous space $ LG /G =\Omega G$ of based loops in
$G$.
%
\item
\labell{2Holes} The moduli space $\M(\Sig^2_0)$ for the two-holed
sphere is diffeomorphic to $LG\times L\g^*$.  Note that these two
factors are of different Sobolev class!  The moment map for the two
$LG$- actions under this identification are $(g,\xi)\mapsto-\xi$ and
$(g,\xi)\mapsto g\cdot\xi$, respectively.
\end{enumerate}
\end{example}
Just as for the case without boundary \cite{at:mo}, the moduli spaces
$\M(\Sig)$ can be described in terms of holonomies.  The fundamental
principle that applies here is that any two flat connections that
define the same holonomy map $\pi_1(\Sig) \ra G$ are related by a
gauge transformation in $\G(\Sig)$.  In the case of a $b$-holed sphere
with $b\ge 1$, i.e.  at least 1 boundary component this leads to the
following description:
\begin{proposition} \labell{HolonomyDescription}\cite{me:lo}
The moduli space $\M(\Sigma_h^{b})$ is equivariantly diffeomorphic to
the smooth submanifold of $G^{2h} \times G^{b}\times (L\g^*)^{b}\ni
(a,c,\xi) $ defined by
$$
 \prod_{i=1}^b \, \Ad_{c_i} \, \Hol(\xi_i)= \prod_{j=1}^{h} \,
 [a_{2j},a_{2j-1}]
$$
where $c_1 = 1$. Here the action of $g= (g_1,\ldots, g_r) \in LG^b$ is
given by
$$ g \cdot a_j = \Ad_{g_1(1)}\, a_j, \ \ g \cdot c_j=g_1(1)\,c_j
g_j(1)^{-1} ,\ \ g \cdot \xi_j=g_j\cdot \xi_j
$$
and the moment map is given by projection to the $(L\g^*)^b$-factor.
\end{proposition}

This description shows for example that the action of 
$(LG)^{b-1}\subset LG^b$ on $\M(\Sigma_h^{b})$ is {\em free}. Setting 
$d_j=\Hol(\xi_j)$ and eliminating $d_1$ one obtains a description 
of the holonomy manifold:

\begin{corollary}
The holonomy manifold $\Hol(\M(\Sigma_h^{b}))$ is $G^{2h+2(b-1)} \ni (a,c,d)$,
with action of $h=(h_1,\ldots,h_b) \in G^{b} $ given by
$$ h \cdot a_j = \Ad_{h_1} \, a_j, \ \ h \cdot c_i=h_1\,c_i\,h_i^{-1},\,\,
   h \cdot d_i=\Ad_{h_j}\,d_i \ \ \ 2 \le i \le b$$
The components of the map $\Hol(\Phi)$ are $ \Hol(\Phi)_j(a,c,d)=d_j$
for $j\ge 2$ and
$$ \Hol(\Phi)_1(a,c,d)= \prod_{j=1}^{h} \, [a_{2j},a_{2j-1}]
 \Big(\prod_{i=2}^b \Ad_{c_i}(d_i)\Big)^{-1}.$$
\end{corollary}
As an application, we have:
\begin{proposition}\labell{EquivariantCohomology}
There is a ring isomorphism in equivariant cohomology (with coefficients in 
$\Z$), $ H^*_{G^{b}}\Hol(\M(\Sigma_h^{b}))=H^*_G(G^{2h+b-1})$
where $G$ acts by conjugation.
\end{proposition}

\begin{proof}
In the description of $\Hol(\M(\Sigma_0^{b}))\cong
G^{2h+2(b-1)}$ given in the Proposition, the action of $G^{b-1}\subset
G^{b}$ is free, and the quotient is just $G^{2h+b-1}$ with $G$ acting
by the adjoint action.  The proposition follows.
\end{proof}

\begin{corollary}
Any $\widehat{LG^{b}}$-equivariant line bundle over $\M(\Sig_h^{b})$
for $b>0$ at level $\levi$ is equivariantly isomorphic to the $k$th
tensor power of the pre-quantum line bundle 
$L(\Sig_h^{b})$.
\end{corollary}

\begin{proof}
It suffices to show that every $\widehat{LG^{b}}$-equivariant line
bundle $L\to \M(\Sig_h^{b})$ at level $0$ is equivariantly trivial.
By Theorem \ref{EquivariantCohomology} and since
$H^1_G(G^{2h+b-1})=H^2_G(G^{2h+b-1})=\{0\}$, it follows that $L/\Om
G^{b}$ is $G$-equivariantly trivial. Consequently, $L$ is
$LG^{b}$-equivariantly trivial.
\end{proof}

\begin{remark}
It is a much deeper fact (see e.g. \cite{at:mo,ku:pg}) that 
in the case $b=0$ every line bundle is a tensor power
of the pre-quantum bundle.
\end{remark}

\subsection{Reduction}
Let $G$ be a compact, connected, simply-connected Lie group and let
$M$ is a Hamiltonian $L(G\times G)$-manifold with proper moment map
$(\Phi_+,\Phi_-)$.  Let $LG$ act on $M$ by the the embedding
$$ \diag: \, LG\to LG\times LG,\,\,g\mapsto (g,I^* g).$$
This action has a moment map at level $0$, $\Phi=\Phi_+ - I^*\Phi_-$.
Notice that $\Phi$ is not proper. In \cite{me:lo} it is shown that if $0$
is a regular value of $\Phi$ the reduced space
$$M\qu \diag(LG):=\Phi^{-1}(0)/ \diag(LG)$$
is a compact orbifold with a naturally induced closed two-form $\om\qu
\diag(LG)$.
One way to see that $M\qu \diag(LG)$ is finite dimensional and compact
is to note that it may also be obtained from the holonomy manifold
$\Hol(M)$.  Indeed, let $\diag:\,G\to G\times G$ be the diagonal
embedding and $\Hol(\Phi)=\Hol(\Phi_+)\cdot\Hol(\Phi_-)$, and define
$\Hol(M)\qu\diag(G)= \Hol(\Phi)^{-1}(e)/\diag(G)$.

\begin{proposition}
Suppose $0$ is a regular value of $\Phi$. The holonomy map
$\Hol:\,M\to \Hol(M)$ descends to a diffeomorphism
$ M\qu \diag(LG) \to \Hol(M)\qu \diag(G)$.
\end{proposition}

\begin{proof}
The condition $m\in\Phinv(0)$ means $\Phi_+(m)=I^*\Phi_-(m)$, hence 
$$\Hol(\Phi_+(m))=\Hol(I^*\Phi_-(m))=\Hol(\Phi_-(m))^{-1}.$$ 
Moreover $(\Om G^2\cdot m)\cap\Phinv(0)=(\diag(\Om G)\cdot m)\cap\Phinv(0)$.
Consequently 
$\Phinv(0)/\diag(\Om G)\cong \Hol(\Phi)^{-1}(e)$ and therefore
$$\Phinv(0)/\diag(\Om G\rtimes G)\cong  \Hol(\Phi)^{-1}(e)/\diag(G).$$ 
\end{proof}

Suppose more generally that $H$ is another simply connected compact,
connected Lie group and that $M$ is a Hamiltonian $L(H\times G\times
G)$-manifold with proper moment map $(\Psi,\Phi_+,\Phi_-)$.  Assume
that $0$ is a regular value of $\Phi= \Phi_+-I^*\Phi_-$ and that the
action of $\diag(LG)$ on the zero level set is free.  Then the reduced
space $M\qu \diag(LG)$ is a Hamiltonian $LH$-Banach manifold with
proper moment map \cite{me:lo}.  (Dropping the freeness assumption
would lead to ``Banach orbifolds.'')

If $L\to M$ is a complex $\widehat{L(H \times G^2)}$-equivariant
line bundle, we can define a reduced bundle
$$ L\qu \diag(LG)=(L|\Phi^{-1}(0))/\diag(LG)$$
which is again at level $k$.  (Note that the diagonal embedding 
$\diag:\,LG\to L(H\times G^2)$ lifts to an embedding into 
$\widehat{L(H\times G^2)}$ so that $\diag(LG)$ acts on $L$.)
If $L$ is a pre-quantum bundle then $L\qu \diag(LG)$ is pre-quantum.

As an example of reduction we have the following result
which we learned from S. Martin (for a proof, see \cite{me:lo}):

\begin{theorem}(Gluing equals reduction)\labell{Factor}
Suppose $\Sig$ is a compact, oriented $2$-manifold that is obtained from a
second $2$-manifold $\widehat{\Sig}$ (possibly disconnected) by gluing
two boundary components $B_\pm\subset\p \widehat{\Sig}$.  Let $LG\hra
LG^b=\G(\p\Sig)$ be the embedding induced by $S^1\hra B_+\times
B_-,\,z\mapsto (z,z^{-1})$. Then $\M(\Sig)$ is the Hamiltonian
quotient
$$ \M(\Sig)= \M(\widehat{\Sig})\qu \diag(LG).$$ 
\end{theorem}

Because every compact oriented 2-manifold $\Sig$ of genus at least 2
admits a pants decomposition, Theorem \ref{Factor} shows that the
moduli space $\M(\Sig)$ can be obtained from products of
$\M(\Sig_0^3)$ by iterated symplectic reductions.

\begin{example} \labell{2Holes2}  It follows from the description
in Example \ref{ex:2Holes} that for any Hamiltonian $LG$-manifold $M$
there is an equivariant symplectomorphism 
$$ M \times \M(\Sig^2_0) \qu \diag(LG) \cong M. $$
\end{example}

As in the finite dimensional case, reductions with respect to
coadjoint orbits are particularly important:

\begin{definition}
Let $(M,\omega,(\Psi,\Phi))$ be a Hamiltonian $LH\times LG$-manifold
with proper moment map. Let $\Alc$ be the alcove for the $G$ and
$\mu\in\Alc$.  If $\mu$ is a regular value of $\Phi$ and the action of
$(LG)_\mu$ on $\Phinv(\mu)$ is free then $0$ is a regular value for
the diagonal action on $M\times LG\cdot(*\mu)$ and the action of
$\diag(LG)$ on the zero level set is also free.  We define the reduced
space at $\mu$ by
$$M_\mu:=(M\times LG\cdot(*\mu))\qu\diag(LG) \cong \Phinv(\mu)/
(LG)_\mu.$$
\end{definition}

For example, the reductions $\M(\Sig_0^b)_{\mu_1,\ldots,\mu_b}$ may
be viewed as moduli spaces of flat connections where the holonomy
around the $j$th boundary component is conjugate to $\Hol(\mu_j)$. 
The holonomy description of these spaces follows from
Proposition \ref{HolonomyDescription}:
\begin{equation}\label{Conreduction}
\M(\Sig_0^b)_{\mu_1,\ldots,\mu_b}= \{(d_1,\ldots,d_b)\in
{\mathcal C}_{* \mu_1}\times\ldots\times {\mathcal C}_{* \mu_b}|\,
\prod_j d_j=e\}/G
\end{equation}
where ${\mathcal C}_{\mu_j}=G\cdot\exp(\mu_j)$ is the conjugacy
class corresponding of $\exp(\mu_j)$ and where the action is the
diagonal action.  For the sake of comparison, note that
$$
\O_{\mu_1}\times \ldots\times \O_{\mu_b}\qu G=
\{(\zeta_1,\ldots,\zeta_b)\in
\O_{\mu_1}\times\ldots\times\O_{\mu_b}|\, \sum_j \zeta_j =0\}/G,
$$
In fact the two spaces are $LG^b$-symplectomorphic,
if the $\mu_i$'s are sufficiently small. This was noticed by 
L. Jeffrey, and also follows from the construction of 
symplectic cross-sections in the following section.

\subsection{Cross-sections} 
\labell{CrossSections}

Hamiltonian $LG$-manifolds with proper moment maps 
behave in many respects like compact Hamiltonian spaces for compact
Lie groups. This is due to the existence of finite dimensional {\em
cross-sections} for the $LG$-action.  The finite collection of
cross-sections behaves like an atlas for a compact, finite dimensional
manifold.
For every open face $\sig$ of the fundamental alcove $\Alc$, we define 
$$ \Alc_\sig:=\bigcup_{\tau \succeq \sig}\tau.$$
Thus $\Alc_\sig$ is the complement in $\Alc$ of the closure of the 
face opposite to $\sig$. 
The flow-out $U_\sig:=(LG)_\sig\cdot\Alc_\sig\subset L\g^*$ is a slice
for the affine $LG$-action;  that is, $LG \cdot U_\sig \cong
LG \times_{(LG)_\sig} U_\sig$. There is a natural $(LG)_\sig$-invariant
decomposition 
$$ TL\g^*|U_\sig=U_\sig\times L\g^*=TU_\sig\oplus U_\sig\times 
(L\g)_\sig^\perp.$$ 
We denote by $(\widehat{LG})_\sig$ the 
restriction to $(LG)_\sig$ of the central extension $\widehat{LG}$. 
\begin{proposition}\labell{Cross}
Let $(M,\omega)$ be a connected Hamiltonian $\LG$-manifold with moment
map $(\Phi,\phi):\,M\to \widehat{L\g}^*=L\g^*\times\R$ such that 
$0\not\in\phi(M)$. Let $\sig\subset\Alc$ be a face of the alcove. 
\begin{enumerate}
\item
The {\em cross-section} $Y_\sig:= (\Phi/\phi)^{-1}(U_\sig)$ is a
smooth Hamiltonian $(\LG)_\sig$-manifold, with the restriction
$(\Phi,\phi)|Y_\sig$ as a moment map. There is a natural $(\LG)_\sig$-invariant
decomposition $TM|Y_\sig=TY_\sig\oplus (L\g)_\sig^\perp$. 
\item
There is a unique $\LG$-invariant closed 2-form and moment map on the
associated bundle $\LG\times_{(\LG)_\sig}Y_\sig$ which restrict to the
given form and moment map on $Y_\sig$. They are obtained by pulling
back $\omega$ resp. $(\Phi,\phi)$ via the embedding
$\LG\times_{(\LG)_\sig}Y_\sig\to \LG\cdot Y_\sig\subset M$.
For each $\tau\prec \sig$, one has canonical Hamiltonian
embeddings
$$
(\LG)_\tau\times_{(\LG)_\sig}Y_\sig\to Y_\tau.
$$
\item
If $(\Phi,\phi)$ is proper as a map into $Lg^* \times \R\backslash\{0\}$
then the cross-section $Y_\sig$ is finite dimensional.
\item
If in addition $\om$ is symplectic, then $Y_\sig$ is a symplectic
submanifold, and the flow-outs $LG \cdot Y_\sig$ are either empty or
connected and dense in $M$.
\end{enumerate}
\end{proposition}

For a proof of these assertions in the case of Hamiltonian
$LG$-manifolds, see \cite{me:lo}.  The extensions to the case of
$\LG$-actions are immediate.  
\begin{remark}
Since the loop group $LG$ is a Hilbert manifold, Proposition 
\ref{Cross} shows that every Hamiltonian $LG$-manifold with proper 
moment map is a Hilbert manifold. 
\end{remark}

One can reconstruct a Hamiltonian $\LG$-manifold $M$ from the
collection of $\{ Y_\sig \}$ and the inclusions $Y_\sig \hra Y_\tau$
for $\tau \prec {\sig}$.  That is, Hamiltonian $\LG$-manifolds with
proper moment maps are equivalent to a collection of finite
dimensional Hamiltonian manifolds with certain compatibility
relations. It is in principle possible to abandon the infinite
dimensional picture and work only with cross-sections; the space $M$
itself serves mainly as a bookkeeping device for the relations between
the cross-sections.  Note in particular that in the case of moduli
spaces of flat connections the cross-sections are independent of the
choice of Sobolev class.

There are some limitations to this point of view - for example a given
invariant almost complex structure on $M$ does not in general preserve the
cross-sections.  Similarly an invariant connection on a given
$\widehat{LG}$-equivariant line bundle is not determined by its
restrictions to the cross-sections.

For the rest of this section we consider only 
Hamiltonian $LG$-manifolds $M$  with proper moment map at level $1$.

\begin{example}
\begin{enumerate}
\item
For any $\mu\in\Alc_\sig$, the cross-section $Y_\sig$ for $M=(LG)\cdot\mu$ 
is the coadjoint orbit for the compact group, ${(LG)}_\sig\cdot\mu$. 
\item
The symplectic cross-sections $Y_\sig$ for the moduli space $\M(\Sig^b_0)$ is 
given in the holonomy picture, 
Proposition \ref{HolonomyDescription}, by imposing the extra 
condition $(\xi_1,\ldots,\xi_b)\in U_\sig$.  These
are the ``twisted extended moduli spaces'' due to 
Jeffrey.
\end{enumerate}
\end{example}

The cross-sections $Y_\sig$ for a Hamiltonian $LG$-manifold 
can also be viewed as submanifolds 
of the holonomy manifold. Indeed, since the intersection $(LG)_\sig\cap \Om G$
is trivial the composition of the horizontal maps 
 in the commutative diagram 
$$ \begin{array}{ccccc}
Y_\sig&\hra& M&\ra&  \Hol(M) \\
\downarrow & & \downarrow & & \downarrow \\
U_\sig & \hra & L\g^*&{\ra}& G 
\end{array} $$
are equivariant embeddings.  Moreover the associated bundles
$$  G\times_{Z_\sig}Y_\sig\to \Hol(M)$$
are open submanifolds. (In fact, one can use these embeddings to define
the manifold structure on $\Hol(M)$.)

Considering the cross-sections as submanifolds of $\Hol(M)$ is useful,
for example, in order to define orientations: Since $G$ and
$Z_\sig\cong (LG)_\sig$ are oriented it follows that any orientation
on $\Hol(M)$ will induce orientations on the cross-sections
$Y_\sig$. Moreover, these orientations are compatible in the sense
that the natural maps $Z_\sig\times_{Z_\tau} Y_\tau\to Y_\sig$ are
orientation preserving. Conversely, compatible orientations on the
cross-sections induce an orientation on $\Hol(M)$.  For example, if
$M$ is a symplectic Hamiltonian $LG$-manifold with proper moment map
$\Phi$ the symplectic orientations on the cross-sections are
compatible.

\subsection{Induction}
Any compact Hamiltonian
$G$-manifold $(M,\omega,\Phi)$ gives rise to a Hamiltonian $LG$-manifold
$\big(\Ind(M),\,\Ind(\omega),\,\Ind(\Phi)\big)$ with proper 
moment map called its {\em induction}.

\begin{proposition}
Let $(M,\omega,\Phi)$ be a Hamiltonian $G$-manifold. There exists 
a unique closed two-form $\Ind(\omega) $ on the associated bundle
\begin{equation}\labell{AssoBundle} 
\Ind(M):=LG\times_G\, M
\end{equation} 
such that the induced $LG$-action on $\Ind(M)$ is Hamiltonian 
with
moment map $\Ind(\Phi)$ at level 1, and such that $\Ind(\Phi)$
and $\Ind(\omega)$ pull back to the given moment map and symplectic
form on $M\subset \Ind(M)$.  
If $M$ is symplectic then the cross-section $Y_{\on{int}(\Alc)}$ is
symplectic. If moreover $\Phi(M)\subset U_{\{0\}}$, 
then $\Ind(M)$ is symplectic.
\end{proposition} 
For a proof, see \cite{me:lo}.  

\begin{example}\labell{InductionCoadjointOrbits}
As before, we denote by $\O_\mu$ the coadjoint $G$-orbit through 
a point $\mu\in\Alc$. If $\mu\in\Alc_{\{0\}}$, the induction 
$\Ind(\O_\mu)$ is just the coadjoint $LG$-orbit through $\mu$. 
For $\mu\not\in\Alc_{\{0\}}$, the induction $\Ind(\O_\mu)$ is not 
symplectic. 
\end{example}

Note that since $LG=\Omega  G\rtimes G$ the associated bundle 
(\ref{AssoBundle}) is simply the trivial bundle $\Omega G\times M$.  
Thus, the holonomy manifold $\Hol(\Ind(M))$ is just $M$ itself, 
and $\Hol(\Ind(\Phi))$ is just the composition $\exp\circ\Phi$.
In particular, it follows from the discussion in section 
\ref{CrossSections} that if $M$ is oriented then all 
cross-sections $Y_\sig$ for $\Ind(M)$ are oriented.
We are mainly interested in the case  
$\sig=\on{int}(\Alc)$ where the cross-sections and their 
orientations are given as follows.
\begin{proposition}
Let $M$ be an oriented Hamiltonian $G$-manifold, with proper 
moment map $\Phi$. For $\sig=\on{int}(\Alc)$, the cross-section 
$Y_\sig$ for the induction $\Ind(M)$ is a disjoint union 
$$  Y_{\on{int}(\Alc)} = \coprod_{w\in\Waff^+} g_w^{-1}\cdot\Phinv(w\cdot\on{int}\Alc) $$
where $g_w\in LG$ is any element representing $w$. 
The orientation on $\Phinv(w\cdot\on{int}\Alc)\subset Y_\sig$ 
is given by $(-1)^{\length(w)}$ times the orientation as a 
subset of the cross-section $\Phinv(\on{int}\,\t^*_+)$ of $M$.
\end{proposition}

\begin{corollary}  \labell{IndRed}
Let $M$ be an oriented Hamiltonian $G$-manifold, with proper moment
map $\Phi$, and $\mu \in \on{int}(\Alc)$.  Then $\mu$ is a regular
value (or not in the image) of $\Ind(\Phi)$ if and only if $w \cdot
\mu$ is a regular value (or not in the image) of $\Phi$, for every $w
\in \Waff^+$.  Furthermore,
$$\Ind(M)_\mu = \coprod_{w \in \Waff^+} (-1)^{\length(w)} M_{w \cdot \mu}$$
\end{corollary}

\begin{proof}  This follows immediately from the proposition since
$$ \Ind(M)_\mu = (Y_\sig)_\mu $$
for $\sig = \on{int}(\Alc)$.
\end{proof}

Instead of proving the proposition in this special case, we will work
out the cross-sections and their orientations in general for arbitrary
$\sig$.  For this we need to consider the set of all images $w\sig$
($w\in\Waff$) which are contained in the positive Weyl chamber $\t_+$.
This set is labeled by the double coset space $W\backslash \Waff /
W_\sig$ where $W_\sig$ is the stabilizer group of $\sig$ in $\Waff$.
The stabilizer group of $w\sig$ in $LG$ is given by
$(LG)_{w\sig}=g\,(LG)_\sig\,g^{-1}$, where $g\in LG$ is any element
with $g\cdot\sig=w\sig$.  For any such $g$, the image
$U_{w\sig}:=g\cdot U_\sig$ is a slice at $w\sig$ for the action of
$LG$. The subset of $\Ind(M)$ given by
$$ Y_\sig^{[w]} := LG_{\sig}\, g^{-1}\cdot\Phinv(U_{w\sig}\cap\g^*)
=g^{-1}\,LG_{w \sig}\cdot\Phinv(U_{w\sig}\cap\g^*)
$$ 
depends only on $w\sig$ and not on the choice of $g$. 

\begin{proposition}
Suppose that $M$ is an oriented Hamiltonian $G$-manifold with equivariant 
map $\Phi:\,M \ra \g^*$. The cross-section 
$Y_\sig$ for the Hamiltonian $LG$-manifold $\Ind(M)$ is a disjoint union 
$$ 
Y_\sig=\coprod_{[w] \in W\backslash \Waff / W_\sig } Y_\sig^{[w]}. 
$$
\end{proposition}

\begin{proof}
We first show the inclusion ``$\supseteq$'': For any $[w] \in
W\backslash \Waff / W_\sig$ and $g\in LG$ such that
$g\cdot\sig=w\sig$, we have
\begin{eqnarray*} 
\Ind(\Phi)(Y_\sig^{[w]})&=&
\Ind(\Phi)(
LG_{\sig}\, g^{-1}\cdot\Phinv(U_{w\sig}\cap\g^*))\\
&\subseteq& (LG)_{\sig}\, g^{-1}\cdot U_{w\sig}= (LG)_\sig\cdot U_\sig=U_\sig.
\end{eqnarray*}
For the opposite inclusion, suppose we are given $x\in Y_\sig$. 
Thus  $x=h\cdot m$ for some $m\in M$, $h\in LG$ with 
$h\cdot \Phi(m)\in U_\sig$.
We may assume with no loss of generality that 
$\mu:=\Phi(m)\in \t_+$. Choose $w\in \Waff^+$ such that 
$\nu:=w^{-1}\mu\in \Alc_\sig$, and choose $g\in LG$ with 
$g\cdot\sig=w\cdot\sig$. 
Then 
$$\mu\in w\Alc_\sig \subset U_{w\sig}\cap\g^*,$$
so that $m\in \Phinv(U_{w\sig}\cap\g^*)$. 
Moreover,
$$g\,h\cdot\mu\in g\cdot U_\sig=U_{w\sig},$$
which implies 
$g\,h\in (LG)_{w\sig}$ since $U_{w\sig}$ is a slice.
\end{proof}

We now consider orientations. Note first that every $w\sig$ is
contained in some fixed open face of the positive Weyl chamber
$\t_+$. 
Therefore all points in $w\sig$ have the same stabilizer group
$G_{w\sig}=G\cap (LG)_{w\sig}$ with respect 
to the coadjoint $G$-action, equal to the 
stabilizer group of that face of $\t_+$. Observe next that
$$U_{w\sig}\cap \g^*=LG_{w\sig}\cdot w\Alc_\sig\cap\g^*=
G_{w\sig}\cdot (w\,W_\sig\Alc_\sig \cap\t_+)$$
is a slice for the coadjoint $G$-action. Hence 
$\Phinv(U_{w\sig}\cap \g^*)\subset M$ 
is an open subset of the cross-section in $M$ corresponding to 
$G_{w\sig}$. We can write $Y_\sig^{[w]}$ as an associated bundle
\begin{equation}\labell{AssociatedBundle} 
Y_\sig^{[w]}=g^{-1}\cdot (LG)_{w\sig}\times_{G_{w\sig}}
\Phinv(U_{w\sig}\cap \g^*).
\end{equation}
Both $\Phinv(U_{w\sig}\cap \g^*)$, being an open subset of 
a cross-section of $M$, and $G_{w\sig}$, being the stabilizer 
of a face of $\t_+$, have canonical orientations. Since 
$g$ is uniquely determined up to right-multiplication by 
an element of $(LG)_\sig$, which is connected, the group 
$(LG)_{w\sig}=g (LG)_\sig g^{-1}$ inherits an orientation from 
the orientation of $LG_\sig$, and finally the associated bundle
(\ref{AssociatedBundle}) obtains a canonical orientation.

\begin{remark}
The orientation of $(LG)_{w\sig}$ depends on $w\sig$.  For example, if
$\sig=\on{int}{\Alc}$ we have $(LG)_{w\sig}=G_{w \sig}=T$ for any
$w\in \Waff$.  The orientation on $(LG)_{w\sig}$
agrees with the given orientation on $T$ (corresponding to the choice
of $\t_+$) if and only if $\length(w)$ is even.
On the other hand the orientation of $G_{w \sig}$ agrees with the
orientation on $T$ for any $w\in\Waff^+$.
\end{remark}

As we explained in the previous section, the cross-sections $Y_\sig$ also
acquire orientations from the orientation on the holonomy manifold 
$\Hol(\Ind(M))\cong M$.

\begin{lemma}\label{Orient1}
The orientations on $Y_\sig$ just described are identical with the 
orientations induced from $\Hol(\Ind(M))\cong M$. 
\end{lemma}

\begin{proof}
To see this, we describe the image $\Hol(Y_\sig)\subset
\Hol(\Ind(M))=M$.  The evaluation map $LG\to G$ gives an embedding
$(LG)_{w\sig}\to G$. The image of this embedding is given by
$Z_{w\sig}:=\on{Ad}_{g(1)}\,Z_\sig$ and contains $G_{w\sig}$.  
By (\ref{AssociatedBundle}), 
$$
\Hol(Y_\sig^{[w]})
= g(1)^{-1}\cdot Z_{w\sig} \times_{G_{w\sig}}\Phi^{-1}(U_{w\sig}
    \cap\g^*).
$$
with an orientation coming from the orientations on 
$Z_{w\sig},\,G_{w\sig}$ and $\Phi^{-1}(U_{w\sig}
    \cap\g^*)$. Since the natural isomorphisms
\begin{eqnarray*}
G\times_{Z_\sig}\Hol(Y_\sig^{[w]})
&\cong&G\times_{Z_{w\sig}}\big(g(1)\cdot \Hol(Y_\sig^{[w]})\big)\\
&\cong&G\times_{Z_{w\sig}} \big(Z_{w\sig}
    \times_{G_{w\sig}}\Phinv(U_{w\sig}\cap \g^*)\big)\\
&\cong&G\times_{G_{w\sig}}\Phinv(U_{w\sig}\cap \g^*)
 \end{eqnarray*}
are orientation preserving, the orientation on $Y_\sig^{[w]}\cong 
\Hol(Y_\sig^{[w]})$ 
just described agrees with the orientation on 
$\Hol(Y_\sig^{[w]})$ induced from the 
embedding $G\times_{Z_\sig}\Hol(Y_\sig^{[w]})\to M$ as an open 
subset.
\end{proof}

\section{Construction of the cobordism}

\subsection{The fusion product of Hamiltonian LG-manifolds}
As we mentioned earlier, Hamiltonian $LG$-manifolds with proper moment
maps behave in many respects like compact Hamiltonian spaces for
compact Lie groups. Now it is well-known that the classical analog
of taking the tensor product of representations is taking the direct
product of Hamiltonian spaces with the diagonal $G$-action. However, 
the direct product $M_1\times M_2$ of two Hamiltonian $LG$-manifolds 
(at level $1$) with proper moment maps has a non-proper moment map 
(at level $2$). There is a different product operation that 
preserves the level and also properness of the moment map.

\begin{definition}
Let $M_1$, $M_2$ be Hamiltonian $LG$-manifolds with proper moment maps
$\Phi_1,\Phi_2$.  Let $\M(\Sig_0^3)$ be the moduli space for the
three-holed sphere. The fusion product $M_1\fus M_2$ is the
Hamiltonian $LG$-manifold obtained as the Hamiltonian quotient
$$ M_1\fus M_2:=\big(M_1\times M_2\times \M(\Sig_0^3)\big)\qu
\diag(LG^2) $$
under the diagonal $LG^2$-action. We denote the resulting moment
map by $\Phi_1\fus \Phi_2:\, M_1\fus M_2\to L\g^*$.
\end{definition}

Observe that $M_1\fus M_2$ is a smooth Banach manifold.  Indeed, since
the $LG^2\subset LG^3$-action on $\M(\Sig_0^3)$ is free, the
corresponding moment map $\M(\Sig_0^3)\to (L\g^*)^2$ is a submersion.
This implies that $0$ is a regular value of the moment map for the
$\diag(LG^2)$-action on $M_1\times M_2\times \M(\Sig_0^3)$, and that
$\diag(LG^2)$ acts freely on the zero level set.  From general
properties of symplectic reduction, it follows that $\Phi_1\fus\Phi_2$
is again proper, and that the fusion product of symplectic Hamiltonian
$LG$-manifolds is symplectic.

In the case of moduli spaces of flat connections we have
$$\M(\Sig_{g_1}^1) \fus \M(\Sig_{g_2}^1) = \M(\Sig_{g_1 + g_2}^1).$$  
This follows immediately from the definition and Theorem
\ref{Factor}.

The category of Hamiltonian $LG$-manifolds with proper moment maps with
product operation $\fus$ may be considered the classical analog to
the fusion ring of $LG$ representations at a given level. 

\begin{proposition}(The classical fusion ring)
\labell{ClassicalFusionRing}
\begin{enumerate}
\item\labell{Part(a)}  For any Hamiltonian $LG$-manifold $M$,
there is an equivariant symplectomorphism $M \fus \Omega G \cong M$.
\item\labell{Part(b)} Let $M_1,\,M_2,M_3$ be Hamiltonian
$LG$-manifolds with proper moment maps.  There are equivariant
symplectomorphisms
$$ M_1\fus M_2 \cong M_2 \fus M_1,$$ 
$$M_1^*\fus M_2^* \cong (M_1\fus M_2)^*.$$
$$ (M_1\fus M_2)\fus M_3 \cong M_1\fus (M_2\fus M_3). $$
\item\labell{Part(c)}
If $LG\cdot\mu$, $LG\cdot \nu$ are coadjoint orbits through
$\mu,\nu\in\Alc$ then $(LG\cdot\mu \fus LG\cdot \nu)_0$ is
a point for $\nu=*\mu$ and empty otherwise.
\end{enumerate} 
\end{proposition}
\begin{proof} 
Assertion (a) follows from
\begin{eqnarray*} 
 M \fus \Omega G &\cong& M \times \M(\Sig^1_0) \times \M(\Sig^3_0) \qu
\diag(LG)^2 \\ &\cong& M \times \M(\Sig^2_0) \qu \diag(LG) \cong M
\end{eqnarray*}
by Example \ref{2Holes2} and Theorem \ref{Factor}.

For any orientation preserving diffeomorphism of $\Sig^3_0$ the
pull-back map on $\A(\Sig)$ induces an equivariant symplectomorphism
of $\M(\Sig^3_0)$.  Therefore, the first assertion in (b) follows from
the existence of an orientation preserving diffeomorphism transposing
two boundary components.  
The second assertion follows from Remark
\ref{Selfadjoint}.  
Associativity follows from Theorem \ref{Factor}, using two 
different ways of cutting a 4-holed sphere into two 
3-holed spheres:
\begin{eqnarray*} 
(M_1\fus M_2)\fus M_3 &\cong& \M(\Sig^3_0) \times
 \M(\Sig^3_0) \times  M_1 \times M_2 \times M_3 \qu \diag(LG)^4 \\
&\cong& \M(\Sig^4_0)  \times  M_1 \times M_2 \times M_3 \qu \diag(LG)^3
\cong
M_1\fus (M_2 \fus M_3).
\end{eqnarray*}

Part (c) follows from 
$$ (LG\cdot\mu \fus LG\cdot \nu)_0\cong \M(\Sig^3_0)_{\mu,\nu,0} \cong
\M(\Sig^2_0)_{\mu,\nu}$$ 
which is a point if $\mu=*\nu$ and empty
otherwise, by the holonomy description of $\M(\Sig^2_0)_{\mu,\nu}$.
\end{proof}

\subsection{Fusion product of holonomy manifolds}
What is the holonomy manifold of a fusion product of Hamiltonian
$LG$-manifolds? Consider the category of isomorphism classes of pairs
$(Q,\Psi)$ where $Q$ is a $G$-manifold and $\Psi:\,Q\to G$ an equivariant
map.  The involution $*$ defines an involution on this
category. Define a fusion product $(Q_1\fus Q_2,\,\Psi_1\fus \Psi_2)$
by setting $Q_1\fus Q_2=Q_1\times Q_2$ with diagonal $G$-action and
$\Psi_1\fus \Psi_2=\Psi_1\cdot\Psi_2$. This fusion product satisfies
axioms analogous to \ref{ClassicalFusionRing}. For example,
commutativity follows by considering the map
$$\tau:\,Q_1\times Q_2\to Q_1\times Q_2,\,\,(q_1,q_2)\mapsto 
(q_1,\Psi_1(q_1)^{-1}\cdot q_2)$$
since $(\Psi_1\fus \Psi_2)\circ \tau = \Psi_2\cdot \Psi_1$. The analog 
to coadjoint orbits $LG\cdot \mu$ ($\mu\in\Alc$) are the conjugacy classes
$\mathcal{C}_\mu=\Ad(G)\cdot\exp(\mu)$. 
\begin{lemma}
Let $(M_i,\om_i,\Phi_i)$ be Hamiltonian $LG$-manifolds with proper moment 
maps. Then 
$$ \Hol(M_1\fus M_2)=\Hol(M_1)\fus \Hol(M_2).$$
\end{lemma}
\begin{proof}  The claim follows from the identities
\begin{eqnarray*} 
\Hol(M_1\fus M_2)&=& \Hol(M_1 \times M_2 \times \M(\Sig^3_0) \qu \diag(LG)^2)
\\ &=& \Hol(M_1)\times \Hol(M_2) \times \Hol(\M(\Sig^3_0)) \qu G^2
\end{eqnarray*}
and the holonomy description given in Proposition
\ref{HolonomyDescription}.
\end{proof}
Note that we also have $\Hol(LG\cdot\mu)=G \cdot \exp(\mu)$ and 
$\Hol(M^*)=\Hol(M)^*$, so that in this sense the holonomy map is a 
homomorphism of classical fusion rings. 
\begin{example}
Let $\mu_1,\ldots,\mu_b\in\Alc$ and $\mathcal{C}_{\mu_j}= G\cdot\exp(\mu_j)$
their conjugacy classes. 
Then (\ref{Conreduction}) may be re-written
$$ 
\M(\Sig_0^b)_{\mu_1,\ldots,\mu_b}=(\mathcal{C}_{\mu_1}
\fus\ldots\fus \mathcal{C}_{\mu_b})\qu G.
$$
\end{example}

\subsection{Induction and fusion product}

The main result of this section is the following relation between
induction and the fusion product:

\begin{theorem}\labell{ClassicalInductionProperty}
Let $(M_i,\omega_i,\Phi_i)$ $(i=1,2)$ be compact Hamiltonian
$G$-manifolds. Then there exists an $LG$-equivariant diffeomorphism
\begin{equation*}
\phi:\, \Ind(M_1)\fus \Ind(M_2) \to \Ind(M_1\times M_2)
\end{equation*}
such that the equivariantly closed 2-forms 
on both sides are cohomologous. That is, there exists
an $LG$-invariant $1$-form $\beta \in 
\Omega^1_{LG}(\Ind(M_1)\fus \Ind(M_2))$ such that 
$$ \Ind(\omega_1)\fus \Ind(\omega_2) = \phi^* \Ind(\omega_1 + \omega_2)
+ \d \beta $$
$$ \Ind(\Phi_1) \fus \Ind(\Phi_2) = \phi^* \Ind(\Phi_1 + \Phi_2)
- \beta^\sharp$$ 
where $\beta^\sharp:\, \Ind(M_1)\fus \Ind(M_2) \to
L\g^*$ is defined by $\l\beta^\sharp,\xi\r =\iota(\xi)\beta$.
\end{theorem}

\begin{proof}
Let $Y_{0,0,0} \subset \M(\Sig_0^3)$ denote the cross-section at
$(0,0,0)\in\Alc\times\Alc\times\Alc$ given in the holonomy description
Proposition \ref{HolonomyDescription} by
\begin{equation}\labell{Holonomy2}
Y_{0,0,0}=\{(c,\eta)\in G^2\times \g^3|\,
\eta_i\in U_{\{0\}} \text{ and }
\exp(\eta_1)\,\Ad_{c_1}\exp(\eta_2)\,\Ad_{c_2}\exp(\eta_3)=1\}.
\end{equation}
Let $\Psi:\,\M(\Sig_0^3)\to (L\g^*)^3$ be the moment map for
$\M(\Sig_0^3)$.  The zero level set is
$\Psi^{-1}(0) \cong G^2$, with two copies of $(LG)_0\cong G$ acting from the 
right and one copy acting diagonally from the left. 

Note that $\Psi^{-1}(0)$ is also the zero-level set for the $G^2\subset G^3$-action 
on the cross-section $Y_{0,0,0}$. 
By the equivariant symplectic normal form theorem, it follows that a neighborhood of
$\Psi^{-1}(0)$ inside $Y_{0,0,0}$ is modeled by a neighborhood of 
the zero section of 
the cotangent bundle $T^*G^2$, 
with two $G$-copies acting by the cotangent lift of the right 
action and one by the  cotangent lift of the left action. More precisely,  
there exists a $G^2$-invariant
neighborhood $V\subset (U_{\{0\}})^2\subset (L\g^*)^2$ of $0$ such that
$\Psi^{-1}(V \times U_{\{0\}})$ is $G^3$-equivariantly
symplectomorphic to a neighborhood $U$ of the zero section in
$T^*G^2$.

In the special case that $\Phi_1(M_1) \times \Phi_2(M_2) \subset V$
then
\begin{eqnarray*} 
\Ind(M_1) \fus \Ind(M_2) &=&  \Ind(M_1) \times \Ind(M_2) \times \M(\Sig^3_0) \qu LG^2 \\
&=&  \Ind(M_1) \times \Ind(M_2) \times \Ind(U)) \qu LG^2 \\
&=& \Ind(M_1 \times M_2 \times U \qu G^2 ) \\
&=& \Ind(M_1 \times M_2) \end{eqnarray*}
as claimed.

To reduce the general case to the case that $\Phi_1(M_1) \times
\Phi_2(M_2) \subset V$ we apply a generalization of the
Duistermaat-Heckman principle. Let $M_j^{(a)}$ denote $M_j$, with
symplectic 2-form multiplied by a factor $a>0$.  By Theorem
\ref{DHApp} below, there exists a family $\varphi_a$ of
$LG$-equivariantly diffeomorphisms $\Ind({M}_1)\fus \Ind({M}_2) \cong
\Ind({M}_1^{(a)})\fus \Ind({M}_2^{(a)})$ such that the equivariant
cohomology class of the equivariant closed 2-form $\Ind(a\om_1) \fus
\Ind(a \om_2)$ varies linearly in $a$.  That is, there exists a closed
$LG$-invariant 2-form $\gamma\in\Om^2(\Ind({M}_1)\fus \Ind({M}_2)$
with moment map $\Psi$ at level $1$, and a family of $LG$-invariant
1-forms $\beta_a\in\Om^1(\Ind({M}_1)\fus \Ind({M}_2))$ such that
$$ \varphi_a^* (\Ind(a\om_1) \fus \Ind(a \om_2)) = \Ind(\om_1) \fus
\Ind(\om_2) +(a-1)\gamma+\d\beta_a,$$
$$ \varphi_a^* (\Ind(a\Phi_1) \fus \Ind(a \Phi_2))=  \Ind(\Phi_1) \fus
\Ind(\Phi_2)+(a-1)\Psi-\beta_a^\sharp$$
where $\langle \beta_a^\sharp,\xi\rangle=\iota(\xi_M)\beta_a$.  

The inductions
$\Ind(M_1^{(a)}\times M_2^{(a)})$ satisfy an analogous property; 
in fact, in this case the 2-form and moment map 
depend linearly on $a$.  For $a$ small enough, $\Phi_1(M_1^{(a)})
\times \Phi_2(M_2^{(a)}) = a\Phi_1(M_1) \times a \Phi_2(M_2) \subset
V$, so that $\Ind({M}_1^{(a)})\fus \Ind({M}_2^{(a)})$ and
$\Ind(M_1^{(a)}\times M_2^{(a)})$ are equivariantly symplectomorphic.
In particular, the slope for the change of cohomology class of the
closed 2-form is the same for both spaces, which proves the theorem.
\end{proof}

An application of Example \ref{Product} (which generalizes immediately
to the $LG$-equivariant setting) shows that $\Ind(M_1)\fus \Ind(M_2)$ and
$\Ind(M_1\times M_2)$ are cobordant as Hamiltonian $LG$-manifolds with
proper moment maps.  (Properness for the moment map $\Phi_N$ for the
cobordism $N$ follows from compactness of $\Hol(N)$.) If $\tau\in\Alc$
is a regular value for $\Phi_N$ this gives a cobordism of compact
orbifolds
$$ (\Ind(M_1)\fus \Ind(M_2))_\tau \sim \Ind(M_1\times M_2)_\tau. $$
If $\tau\in\on{int}\Alc$ the reduction on the right hand side is given
by Corollary \ref{IndRed}, and we obtain an oriented cobordism
\begin{equation} \labell{QuotCob}
 (\Ind(M_1)\fus \Ind(M_2))_\tau \sim \coprod_{w\in \Waff^+}
(-1)^{\length(w)} (M_1\times M_2)_{w \tau}.
\end{equation}
Thinking of the cobordism in (\ref{QuotCob}) as a quotient of the
finite-dimensional cross-section $\Phinv_N(\on{int}(\Alc))$, the
perturbation argument of Ginzburg-Guillemin-Karshon \cite{gi:co} 
shows that it suffices to assume that $\tau$ 
is a regular value for moment maps of the $LG$-actions on the ends of the 
cobordism.

We are particularly interested in the case where $M_1=\O_{*\mu}$ and
$M_2=\O_{*\nu}$ are coadjoint orbits through
$*\mu,*\nu\in\Alc_{\{0\}}$:  
\begin{theorem}\labell{ClassicalCobordism}
Let $\mu,\nu\in\Alc_{\{0\}}$ and $\tau\in\on{int}(\Alc)$ such that for
every $w \in \Waff^+$, $w \cdot \tau$ is a regular value (or not in
the image) of the moment map for the diagonal action of $G$ on
$\O_{*\mu}\times\O_{*\nu}$, and $\tau$ is a regular value (or not in
the image) of the moment map for the action of $LG$ on
$\M(\Sig^3_0)_{\mu,\nu}$.  Then there is an oriented orbifold
cobordism
\begin{equation}\labell{Cob}
\M(\Sig_0^3)_{\mu,\nu,\tau}\sim \coprod_{w\in \Waff^+}
(-1)^{\length(w)}(\O_{*\mu}\times\O_{*\nu} \times \O_{w(*\tau)} ) \qu G
\end{equation}
where the symplectic forms extend to a closed two-form on the
cobordism. In the case of $G=\on{SU}(n)$, the reduced spaces 
on both sides are smooth manifolds and the cobordism is a 
cobordism of smooth manifolds. 
\end{theorem}
\begin{proof}
The induced spaces $\Ind(\O_{*\mu})$ and
$\Ind(\O_{*\nu})$ are simply the coadjoint $LG$-orbits through
$*\mu,\,*\nu$, hence $(\Ind(\O_{*\mu})\fus \Ind(\O_{*\nu}))_\tau$ is the
moduli space $\M(\Sig_0^3)_{\mu,\nu,\tau}$.  Therefore the cobordism
is a special case of (\ref{QuotCob}).

In the case  $G=\on{SU}(n)$, 
every discrete stabilizer for the $LG$-action on $\M(\Sig^3_0)_{\mu,\nu}$
is equal to the center $Z(G)$. This follows from the holonomy description, 
as explained in \cite{me:wi}, Remark 4.3. Consequently the action 
of $T/Z(G)$ on $\Phi_N^{-1}(\tau)$ is free in this case and $N_\tau$ is 
smooth.
\end{proof}

The cobordism given in this Theorem can be carried further.  
Consider $\g/\t$ as a Hamiltonian $T$-manifold, with symplectic
structure corresponding to the $T$-invariant Hermitian 
structure defined by the metric and the choice of $\t_+$. 
By \cite{me:wi}, Section 2
there are cobordisms
$$ (\O_{*\mu} \times \O_{*\nu})_{w\tau}\sim \coprod_{w_1,\,w_2\in W}
(-1)^{\length(w_1 w_2)} (\g/\t)_{-w_1*\mu-w_2*\nu+w\tau}.$$
Combining this with (\ref{Cob}) gives: 
\begin{theorem}\labell{ToricVarieties}
Let $\mu,\nu\in\Alc_{\{0\}}$ and $\tau\in\on{int}(\Alc)$ such that
$\tau$ is a regular value (or not in the image) of the action of $LG$
on $\M(\Sig^3_0)_{\mu,\nu}$ and for all $w \in \Waff^+$ and
$w_1,w_2\in W$, $-w_1*\mu-w_2*\nu+w\tau$ is a regular value (or not in
the image of) the $T$-action on $\g/\t$. Then
$\M(\Sig_0^3)_{\mu,\nu,\tau}$ is cobordant to a disjoint union of
toric varieties
$$\M(\Sig_0^3)_{\mu,\nu,\tau}\sim 
 \coprod_{w\in \Waff^+,\,w_1,\,w_2\in W}
(-1)^{\length(w\,w_1\, w_2)} (\g/\t)_{-w_1*\mu-w_2*\nu+w\tau}.
$$
\end{theorem}

\begin{remark}
Theorem \ref{ClassicalCobordism} and the following discussion
generalizes easily to the case of more than 3 markings. Given
$\mu_1,\ldots,\mu_{b-1}\in\Alc_{\{0\}},\mu_b\in\on{int}(\Alc)$ we have
$$ \M(\Sig_0^b)_{\mu_1,\ldots,\mu_{b}}\sim \coprod_{w\in
\Waff^+}
(-1)^{\length(w)} ( \O_{*\mu_1}\times\ldots\times
\O_{*\mu_{b-1}}\times\O_{w(*\mu_{b})}) \qu G.
$$
\end{remark}

\section{Applications}
\subsection{Mixed Pontrjagin numbers for $\M(\Sig_0^3)_{\mu_1,\mu_2,\mu_3}$}
In this section we use Proposition \ref{ToricVarieties} to compute the
mixed Pontrjagin numbers for the moduli space of the three-holed
sphere. Recall that for regular
$\mu=(\mu_1,\mu_2,\mu_3)\in\on{int}\Alc^3$, the moduli space $
\M(\Sig_0^3)_{\mu_1,\mu_2,\mu_3}$ has complex dimension equal to 
$k=\f{1}{2}(\dim G-3\dim T)$.  Given an invariant polynomial
$$p\in S(\mathfrak{u}(k)^*)^{\on{U}(k)}\cong
\C[x_1,\ldots,x_k]^{S_k}$$
we can define the corresponding Chern class
$$ p(\M(\Sig_0^3)_{\mu_1,\mu_2,\mu_3})\in
   H^*(\M(\Sig_0^3)_{\mu_1,\mu_2,\mu_3})$$
and the mixed Chern number
$$
\int_{\M(\Sig_0^3)_{\mu_1,\mu_2,\mu_3}}p(\M(\Sig_0^3)_{\mu_1,\mu_2,\mu_3})
\exp(\om_{\mu_1,\mu_2,\mu_3})$$
where $\om_{\mu_1,\mu_2,\mu_3}$ is the symplectic form on
$\M(\Sig_0^3)_{\mu_1,\mu_2,\mu_3}$.\\  
The Pontrjagin ring of $\M\big(\Sig_0^3)_{\mu_1,\mu_2,\mu_3}\big)$ is the subring
generated by polynomials
$$p\in S(\mathfrak{o}(2k)^*)^{\on{O}(2k)}\cong \C[x_1^2,\ldots,x_k^2]^{S_k}
\subset \C[x_1,\ldots,x_k]^{S_k}
$$
The {\em mixed Pontrjagin numbers} corresponding to these polynomials
are invariants under oriented cobordism of orbifolds with closed
2-forms.  Proposition \ref{ToricVarieties} reduces the computation 
of the Pontrjagin
numbers of $\M(\Sig_0^3)_{\mu_1,\mu_2,\mu_3}$ to that for the reduced
spaces of $\g/\t$, which are toric varieties.  
To any symmetric polynomial $p\in \C[x_1,\ldots,x_k]^{S_k}$ we associate
a polynomial on $\t$ as follows. Let $n=\f{1}{2}\dim(G/T)$. By 
the canonical extension map 
$$ \C[x_1,\ldots,x_k]^{S_k}\to \C[x_1,\ldots,x_n]^{S_n}$$
(sending the $l$th elementary symmetric polynomial in $k$ variables to that in
$n$ variables) we can view $p$ as a polynomial on
$\R^n$. Choose an ordering $\alpha_1,\ldots,\alpha_n$ of the 
positive roots of $G$ to identify $\oplus_{\alpha\in\mathfrak{R}_+}
\R\cong\R^n$ and $\g/\t\cong\C^n$. 
The action of $T$ on $\g/\t$ gives a map
$T\to \on{U}(1)^n$.  We define the polynomial $\phi^*\,p\in S(\t^*)$
to be the pull-back of $p$ by the tangent map $\phi:\, \t \ra
\u(1)^n=\R^n$. The polynomial $\phi^*\,p$ defines a  
constant coefficient differential operator $(\phi^*\,p)(\f{\p}{\p\mu})$ 
on $\t^*$.

Let $\kappa$ be the push-forward of the characteristic measure on the
positive orthant $\R^n_+$ under the map $\phi^*:\,\R^n\to \t^*$ 
dual to $\phi$. 
By comparing with the given Lebesgue measure on
$\t^*$ we can consider $\kappa$ as a (piecewise polynomial) function
on $\t^*$.
\begin{lemma}\label{GmodT}
The symplectic volume of any reduced space $(\g/\t)_\mu$
at a regular value $\mu\in\t^*$ is given by
$$ \Vol\big((\g/\t)_\mu\big) = \int_{(\g/\t)_\mu}\,\exp(\om_\mu)=
\f{\# Z(G)}{\Vol(T)} \kappa(-\mu).$$
The mixed Chern
number corresponding to an invariant polynomial $p$ is given by
application of the differential operator $\phi^*p(\f{\p}{\p \mu})$ to
the volume function:
\begin{equation}\label{Chern}
 \int_{(\g/\t)_\mu}\,p\big((\g/\t)_\mu\big)\,\exp(\om_\mu)
=(\phi^*p)(\f{\p}{\p \mu})\Vol\big((\g/\t)_\mu\big).
\end{equation}
\end{lemma}

\begin{proof}
The proof of this result can be found e.g. in \cite{gu:mo} or 
\cite{gu:dh}. We recall the argument for convenience of 
the reader.

The Duistermaat-Heckman measure for the standard $U(1)^{n}$-action 
on $\g/\t$ is given by the characteristic function of the negative orthant 
$-\R_+^{n}\subset \R^{n}$. 
Therefore the Duistermaat-Heckman measure 
for the $T$-action is $\mu\mapsto \kappa(-\mu)$, which proves 
the volume formula. The factor $\# Z(G)$ is due to the
fact that the generic stabilizer for the $T$-action on $\g/\t$ is
equal to $Z(G)$.

To prove (\ref{Chern}) let $\Phi:\ \g/\t \ra \t^*$ denote the
moment map for the $T$-action and $\pi:\, \Phinv(\mu) \ra
(\g/\t)_\mu$ and $\iota:\,\Phinv(\mu)  \ra \g/\t$ the projection and
inclusion. 
Then $ \pi^* T(\g/\t)_\mu \oplus \t_\C = \iota^* \, \g/\t$
and therefore
$$
T(\g/\t)_\mu\oplus \t_\C=
\bigoplus_{\alpha\in\mathfrak{R}_+} L_\alpha $$
where $L_\alpha = Z \times_T \C_\alpha$ (see e.g. \cite[p. 59]{gu:mo}).
The Chern classes of $(\g/\t)_\mu$ are symmetric polynomials in the
first Chern classes $c_\alpha=c_1(L_\alpha)$.  
The classes $c_\alpha$ are related to the first Chern class $c \in
\Omega^2\big((\g/\t)_\mu,\t\big)$ of the torus bundle $Z \ra (\g/\t)_\mu$ by $
(c_{\alpha_1}, \ldots, c_{\alpha_n}) = \phi \circ c$.  Thus
$ (\phi^*p)(c) = p(c_{\alpha_1}, \ldots, c_{\alpha_n})
=p\big((\g/\t)_\mu\big).$ This proves the Theorem since by the Duistermaat-Heckman
Theorem,
$$ \phi^*p\,(\f{\p}{\p \mu})\Vol\big((\g/\t)_\mu\big)=\int_{(\g/\t)_\mu}
\phi^*p(c)\exp(\om_\mu).$$
\end{proof}

\begin{theorem} \labell{VolPontr}
For $\mu=(\mu_1,\mu_2,\mu_3)\in\on{int}
\big(\Alc^3\cap \Phi(\M(\Sig_0^3))\big)$, the volume of the
moduli space $ \M(\Sig_0^3)_{\mu_1,\mu_2,\mu_3}$ is given by
\begin{equation}\label{SteinbergVol}
(-1)^{\f{1}{2}\dim G/T} \f{\# Z(G)}{\Vol(T)}\,
\sum_{l\in\Lambda}\sum_{w_1,\,w_2\in W}
(-1)^{\length(w_1\,w_2)}\kappa(w_1\mu_1+w_2\mu_2+\mu_3+l).
\end{equation}
If in addition
$\mu$ is a regular value for the moment map then the
mixed Pontrjagin number corresponding to an invariant polynomial $p\in
S^*(\mathfrak{o}(2k))^{\on{O}(2k)}$ (where $k=\f{1}{2}(\dim G-3\dim
T)$) is given by application of the differential operator $\phi^*
p(\f{\p}{\p \mu_3})$ to the volume function:
$$
\int_{\M(\Sig_0^3)_{\mu_1,\mu_2,\mu_3}}p(\M(\Sig_0^3)_{\mu_1,\mu_2,\mu_3})
\exp(\om_{\mu_1,\mu_2,\mu_3})
=(\phi^*p)(\f{\p}{\p \mu_3})\Vol(\M(\Sig_0^3)_{\mu_1,\mu_2,\mu_3}) 
$$
\end{theorem}

\begin{proof}
Using Proposition \ref{ToricVarieties}, Lemma \ref{GmodT}
and the invariance of mixed 
Pontrjagin numbers under oriented cobordism,  the mixed Pontrjagin 
numbers are given by applying $(\phi^*p)(\f{\p}{\p \mu_3})$ to 
$$ \Vol(\M(\Sig_0^3)_{\mu_1,\mu_2,\mu_3 }) = \f{\#
Z(G)}{\Vol(T)}\sum_{w\in \Waff^+,\,w_1,\,w_2\in W}
(-1)^{\length(w\,w_1\, w_2)}\kappa(w_1\mu_1+w_2\mu_2-w*\mu_3)
$$
using the $*$-invariance of $\kappa$. 
The sum over 
$\Waff^+ =\Waff/W\cong \Lambda$ can be replaced 
by a sum over the integral lattice
$\Lambda\subset \Waff$, using the $W$-invariance of the function
$$\nu\mapsto \sum_{w\in W}(-1)^{\length(w)} \kappa(w\mu-\nu).$$
Using $- * \mu = w_0 \mu$ where $\length(w_0) = \dim(G/T)$ gives 
(\ref{SteinbergVol}). 
By continuity the volume formula also holds at non-regular values 
$\mu\in\on{int}\big(\Alc^3\cap\Phi(\M(\Sig_0^3))\big)$.  
\end{proof}

\subsection{Witten's volume formulas}

In this section we briefly outline how to obtain Witten's formulas for
the symplectic volumes of the moduli spaces.  Further details can be
found in \cite{me:lo}.
We label the irreducible $G$-representations by their dominant weights
$\lambda\in\Lambda^*_+:=\Lambda^*\cap \t_+$ and let
$\chi_\lambda:\,G\to \C$ denote the character and
$d_\lambda=\chi_\lambda(e)$ the dimension.  Let 
$\Vol(G)$ be the Riemannian volumes of $G$ with respect to the
normalized inner product on $\g$. For any $\mu\in\Alc$ let 
$\mathcal{C}_\mu\subset G$ the conjugacy class parametrized of 
$\exp(\mu)$ 
and $\Vol(\mathcal{C}_\mu)$ its Riemannian volume.
\begin{theorem}[Witten Formula] \labell{WitThm}
Suppose $2h+b\ge 3$. Let $\mu=(\mu_1 ,\ldots,\mu_b)\in\Alc^b$ be such
that the level set $\Phinv(\mu)$ contains a connection with stabilizer
$Z(G)$.  Then the volume $
\Vol\big(\M(\Sig^b_h,\mu_1,\ldots,\mu_b)\big)$ of the moduli space of
the 2-manifold $\Sig^b_h$ with fixed holonomies $\mu_1,\ldots,\mu_b$
is given by the formula
\begin{equation} \labell{WitForm}
\# Z(G)\, \Vol(G)^{2h-2} \prod_{j=1}^b
\Big(\Vol(\mathcal{C}_{\mu_j})\!\!
\prod_{\stackrel{\alpha\in\mathfrak{R}^{+}}{\l\alpha,\mu_j\r\not\in\{0,1\}
}} \!\!2\sin (\pi \l\alpha,\mu_j\r) \Big)
\sum_{\lambda\in\Lambda^*_+}\f{1}{d_\lambda^{2h-2+b}} \prod_{j=1}^b
\chi_\lambda(e^{\mu_j}).  \end{equation}
In particular, if $\Sig_h^0$ has no boundary:
$$ 
\Vol(\M(\Sig_h^0))=\#\,Z(G)
\Vol(G)^{2h-2}
\sum_{\lambda\in\Lambda^*_+}
\f{1}{d_\lambda^{2h-2}}.
$$
\end{theorem}
We sketch the main ideas, referring to \cite{me:lo} for details.
Witten's formula for $\M(\Sig_0^3)_{\mu_1,\mu_2,\mu_3}$ with
$\mu=(\mu_1,\mu_2,\mu_3) \in \on{int}(\Alc^3\cap\Phi(\M(\Sig_0^3))$
\begin{equation}\labell{Witten1}
\Vol(\M(\Sig_0^3)_{\mu_1,\mu_2,\mu_3})= \# Z(G) \f{\Vol(G)}{\Vol(T)^3}
\sum_{\lambda\in\Lambda^*_+} \f{1}{d_\lambda}\prod_{j=1}^3
\chi_\lambda(e^{\mu_j}) \prod_{\alpha\in\mathfrak{R}^+}\,2\,
\sin(\pi\l\alpha,\xi\r).
\end{equation}
can be proved by applying the the Poisson summation and Weyl character
formulas to (\ref{SteinbergVol}).  Formulas for the general case are
obtained by gluing.  Let $\Sig=\Sig_h^b$ be obtained from a possibly
disconnected 2-manifold $\hat{\Sig}$ by gluing two boundary components
$B_\pm\subset\p\Sig$, and $\mu_1,\ldots,\mu_b\in \on{int}\,\Alc$ such
that $\M(\Sig,\mu_1,\ldots,\mu_b))$ contains at least one connection
with stabilizer $Z(G)$.  Then
$$ \Vol(\M(\Sig,\mu_1,\ldots,\mu_b))= \f{1}{k}\int_{\Alc}
\Vol(\M(\hat{\Sig},\mu_1,\ldots,\mu_b,\nu,*\nu))|\d\nu|
$$
(see e.g.  Jeffrey-Weitsman \cite{jw:va}).  Here the measure $|\d\nu|$
on $\Alc\subset \t$ is the normalized measure for which $\t/\Lambda^*$
has measure 1, and $k=1$ if $\hat{\Sig}$ is connected and equal to $\#
Z(G)$ if $\hat{\Sig}$ is disconnected.
Choosing a pants decomposition for $\Sig$ and using the orthogonality
relations of the characters $\chi_\lam$ and the Weyl integration
formula, carrying out the integrations gives the formula
(\ref{WitForm}) for $\mu=(\mu_1,\ldots,\mu_b) \in\on{int}(\Alc)^b$.  The
volumes for arbitrary $\mu$ are computed by studying the
limit as $\mu$ approaches the boundary of $\Alc^b$.\\

The formula (\ref{WitForm}) was proved in most cases by Witten
\cite{wi:qg}.  Alternative proofs and extensions were given in
\cite{li:hk,li:h2,je:in,jw:va}.

\section{A Duistermaat-Heckman principle}

In this section we prove the Duistermaat-Heckman result 
used in the proof of Theorem \ref{ClassicalInductionProperty}.
We follow the strategy from Section \ref{CrossSections}, carrying out
all constructions in finite-dimensional cross-sections.  As a first
step we construct invariant tubular neighborhoods which are suitably adapted 
to the cross-sections.
Let $(M,\om)$ be a Hamiltonian $\LG$-manifold, with proper moment map
$(\Phi,\phi):\,M\to L\g^*\times\R\backslash\{0\}$.  Let $\sig$ be a
face of $\Alc$ and $Y_\sig$ the corresponding cross-section. Recall
from Proposition \ref{Cross} that there is a canonical isomorphism
\begin{equation} \label{Splitting}
 TM | Y_\sig \cong TY_\sig \oplus \ (L\g)_\sig ^\perp.
\end{equation}
\begin{definition} 
Let $A\subset Y_\sig$ be a subset. 
\begin{enumerate}
\item
An $\LG$-invariant metric $g$ on $M$ is called adapted to $Y_\sig$ over $A$ 
if
there is an open neighborhood $O\subset Y_\sig$ of $A$ such that the
splitting \eqref{Splitting} is orthogonal over $TM|O$, and the induced 
metric on the trivial bundle $O\times (L\g)_\sig^\perp$ agrees with
the $L^2$-metric on  $  (L\g)_\sig^\perp \subset L\g$. 
\item
An $\LG$-invariant  differential form $\alpha$ on $M$ is called adapted to
$Y_\sig$ over $A$ if there is an open neighborhood $O\subset Y_\sig$
of $A$ such that for all $x\in O$, the  subspace 
$(L\g)_\sig ^\perp \subset 
T_xM$ is contained in the kernel of $\alpha$.
\end{enumerate}
\end{definition}
Observe that
adapted metrics or differential forms can be reconstructed over 
$LG\cdot O$ from their restrictions to $O$. If a metric or differential 
form is adapted to $Y_\sig$ over $A\subset Y_\sig$ and if $\tau\preceq\sig$, 
then 
it is also adapted to $Y_\tau$ over $(\LG)_\tau\cdot A\subset Y_\tau$.

It is in general 
impossible in general to choose a metric $g$ on $M$ which is adapted 
to all cross-sections $Y_\sig$ over all of $Y_\sig$.
For this it is necessary to choose smaller subsets
$A\subset Y_\sig$.
We first choose a subdivision $\Alc=\bigcup_\sig \Alc_\sig'$ as
indicated in Figure \ref{AlcoveFig}. The polytopes $\Alc_\sig'$ are
given explicitly as follows. For any face $\sig$ of $\Alc$ let
$C_\sig\subset\t$ be the cone generated by outward-pointing normal
vectors to facets containing $\sig$.  Choose $\eps\in (0,1)$ and $\mu
\in \on{int}(\Alc)$, and define $\Alc_\sig'$ as the intersection
$$\Alc_\sig'=\Alc\cap \big(\eps \mu + 
(1-\eps)(\ol{\sig}-\mu)+ C_\sig\big).$$
\begin{figure}[htb] 
\begin{center}
\setlength{\unitlength}{0.00033333in}
\begingroup\makeatletter\ifx\SetFigFont\undefined%
\gdef\SetFigFont#1#2#3#4#5{%
  \reset@font\fontsize{#1}{#2pt}%
  \fontfamily{#3}\fontseries{#4}\fontshape{#5}%
  \selectfont}%
\fi\endgroup%
{
\begin{picture}(4824,5439)(0,-10)
\path(612,5262)(612,462)(4812,2862) (612,5262)(612,5262)
	\path(912,1062)(1060,720) \path(912,1062)(612,1062)
	\path(912,4662)(612,4662) \path(912,4662)(1100,4950)
	\path(3987,2862)(4172,3172) \path(3987,2862)(4152,2520)
	\path(912,4662)(912,1062)(3987,2862) (912,4662)(912,4662)
	\path(912,4662)(912,1062)(3987,2862)
	(912,4662)(912,4662)
\end{picture}
}
\caption{Subdivision of the fundamental alcove for $G=SU(3)$
\label{AlcoveFig}}
\end{center}
\end{figure}
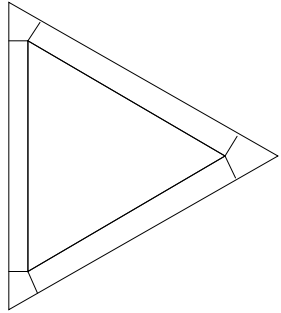\\
Set 
\begin{equation} \labell{ResCrossSec}
Y_\sig' := (\Phi/\phi)^{-1}(\Alc_\sig').\end{equation}
We will construct a metric $g$ on $M$ which is adapted to every 
$Y_\sig$ over the subset $\bigcup_{\tau\succeq\sig}Y_\tau'$. 
\begin{lemma}
\labell{Adapted}
There exists an $G$-invariant partition of unity $\{\rho_\sig\}$
on $G$, subordinate to the cover $G\cdot\Alc_\sig$ such that for
each face $\sig$,\,
\begin{equation} \labell{sum}
 \sum_{\,\tau\succeq \sig
}\rho_\tau=1 \end{equation}
on an open neighborhood of $G \cdot \bigcup_{\,\tau\succeq
\sig} \Alc_\tau'$.
\end{lemma}

\begin{proof} 
Observe first that the conditions \eqref{sum} imply that for every $\sig$, 
\begin{equation} \labell{support}
\on{supp}(\rho_\sig)\subset\bigcup_{\tau\preceq\sig}
G \cdot \Alc_\tau'.
\end{equation}
The proof of Lemma \ref{Adapted} is
by induction, starting from $\sig=\on{int}(\Alc)$.
Suppose by induction that we have constructed a collection $\{
\rho_\sig\}_{\dim\sig > k}$ of non-negative invariant functions with
$\on{supp}(\rho_\sig)\subset G\cdot \Alc_\sig$ such that 
$\sum_{\dim\sig>k} \rho_\sig\le 1$ and 
\begin{enumerate}
\item 
$\on{supp}(\rho_\sig)\subset \bigcup_{\tau\preceq\sig}
G \cdot \Alc_\tau'$ and
\item $\sum_{\,\tau\succeq \sig }\rho_\tau=1$ 
on an open neighborhood of
$G \cdot \bigcup_{\,\tau\succeq \sig} \Alc_\tau'$.
\end{enumerate}
Let $V_k$ be an invariant open neighborhood 
of $G\cdot \bigcup_{\dim\tau>k}\Alc_\tau'$ with
$\big(\sum_{\dim\tau>k}\rho_\tau\big)|V_k=1$. If $\sig_1,\sig_2$ are 
disjoint faces of dimension $k$, we have by construction
$$ (G\cdot\Alc_{\sig_1}'\backslash V_k)\cap
(G\cdot\Alc_{\sig_2}'\backslash V_k)=\emptyset.$$
Therefore we can choose non-negative 
invariant functions $f_\sig\le 1$ for every $k$-dimensional face $\sig$ 
with 
$$ \on{supp}(f_\sig)\subset \bigcup_{\tau \preceq \sig} G\cdot \Alc_\tau',$$
such that the supports of the various $f_\sig$'s are disjoint and such
that $f_\sig=1$ on a neighborhood of $ G\cdot\Alc_\sig'\backslash V_k$.
Then $\rho_\sig=f_\sig\,(1-\sum_{{\tau}\succ \sigma}
\rho_\tau)$
satisfies Condition a and b for all faces $\sig$ with $\dim\sig\ge k$, 
and $\sum_{\dim\tau \ge k}\rho_\tau\le 1$.
\end{proof}

\begin{lemma} \label{metric}
There exists an $\LG$-invariant Riemannian metric on $M$ such
that for each face $\sig$ of $\Alc$ the canonical splitting
(\ref{Splitting}) is adapted to $Y_\sig$ over $\bigcup_{{\tau} 
\succeq \sig} Y_\tau'$. 
\end{lemma}

\begin{proof}
For each $\sig$ choose an $(\LG)_\sig$-invariant Riemannian metric on
$Y_\sig$ and let $g_\sig$ be the corresponding adapted metric on $LG
\cdot Y_\sig$. Let  $\{ \rho_\sig \}$ be a partition of unity as in Lemma
\ref{Adapted}. Then 
$$g=\sum_{\sig} \big((\Hol\circ\Phi)^*\rho_\sig\big) \,g_\sig$$ 
has the required property. 
\end{proof}

Using adapted metrics we can prove the tubular neighborhood 
theorem.

\begin{lemma} \labell{TubNbd}
Let $(X,\om)$ be a Hamiltonian $\LG$-manifold, 
with proper moment map
$(\Phi,\phi):\, X \ra Lg^*\times \R\backslash\{0\}$, 
and suppose $s \in \R\backslash\{0\} $ is a regular
value of $\phi$.  Then for $\eps$ sufficiently small there exists an
equivariant diffeomorphism
\begin{equation}\label{tubular}
F:\,\phinv(s) \times (s-\eps,s+\eps)\to 
 \phinv\big((s-\eps,s+\eps)\big)  
\end{equation}
such that $(F^*\phi)(x,t)=t$, and $F$ is adapted to the cross-sections
in the sense that for each face $\sig \subset \Alc$, there is an open
neighborhood $V_\sig \subset Y_\sig \cap\phinv(s)$ of $Y_\sig' \cap \phinv(s)$
with
\begin{equation}\label{tu}
F\big(V_\sig \times (s-\eps,s+\eps)\big) \subset Y_\sig .
\end{equation}
\end{lemma}

\begin{proof} 
Let $g$ be a Riemannian metric as in Lemma \ref{metric}. 
For $\eps$ sufficiently small and for sufficiently small open
neighborhoods $V_\sig \subset Y_\sig \cap \phinv(s)$ of $Y_\sig' \cap
\phinv(s)$, geodesic flow on $Y_\sig$ in the normal direction to
$Y_\sig \cap \phinv(s)$ defines an $(\LG)_\sig$-equivariant diffeomorphism
$$
\bar{F}_\sig:\, \,V_\sig \times (s-\eps,s+\eps)\hra Y_\sig.
$$
Since $\f{\p}{\p t}\bar{F}_\sig^*\phi(x,t)>0$ there is a unique
reparametrization $\varphi_\sig(x,t)$, defined for $t$ sufficiently
close to $s$, such that $F_\sig(x,t)=\bar{F}_\sig(x,\varphi_\sig(x,t))$
satisfies $F_\sig^*\phi(x,t)=t$.  By construction $F_\sig$ is
$(\LG)_\sig$-equivariant and for $\eps$ sufficiently small extends to
an $\LG$-equivariant diffeomorphism
$$ F_\sig: \ \LG \cdot V_\sig \times
(s-\eps,s+\eps) \hra \LG\cdot Y_\sig .$$
Choosing $\eps$ and $V_\sig$ smaller if necessary we can assume that
for every $\sig$, the metric $g$ is adapted to $Y_\sig$ over the image
of ${F}_\sig$, hence over
$\bigcup_{\tau\succeq\sig}\on{im}({F}_\tau)$. This means that for
$\tau\preceq\sig$, the image $\on{im}(F_\sig)$ is a totally geodesic
submanifold of $(\LG)_\tau\cdot\on{im}(F_\sig)\subset Y_\tau$, 
and the above geodesic flows and reparametrizations coincide. 
Consequently the maps $F_\sig$ patch together to a diffeomorphism $F$
with the required properties.
\end{proof}

There are a number of forms of the Duistermaat-Heckman principle which
hold in the $LG$-equivariant setting.  Here we need only the following
special case.

\begin{proposition} 
\labell{DH} Let $(X,\om,(\Phi,\phi))$ be a Hamiltonian 
$\widehat{LG}$-manifold as in Lemma \ref{TubNbd}. Suppose that 
$s$ is a regular value for $\phi$ and that the central $\UU$ acts freely 
on $\phi^{-1}(s)$.  
Let $F$ and $\eps$ be as in Lemma \ref{TubNbd}. For $t\in (s-\eps,s+\eps)$ 
let $(X_t,\om_t,\Phi_t)$ the reduced space for the central $\UU$ and 
$f_t:\,X_s\to X_t$ the diffeomorphisms induced by $F(\cdot,t)$. 
 Then $f_t^*{\om}_t+2\pi i f_t^*{\Phi}_t$ 
varies linearly in cohomology, with slope $\gamma+2\pi i\Psi$. 
That is, there exists an $LG$-invariant closed 2-form $\gamma\in\Om^2(X_s)$
with moment map at level $s$, $\Psi:\,M\to L\g^*$, and 
a family of $LG$-invariant 1-forms $\beta_t\in\Om^1(X_s)$ such that
$$ f_t^*{\om}_t=\om_s+(t-s)\gamma+\d\beta_t,$$
$$ f_t^*{\Phi}_t={\Phi}_s+(t-s)\Psi-\beta_t^\sharp$$
The form $\gamma$ is the curvature of an $LG$-invariant connection on 
$\phinv(s)$. 
\end{proposition}

\begin{proof}  
For $t\in (s-\eps,s+\eps)$ let $F_t=F(t,\cdot):\ \phinv(s) \to M$
be the family of equivariant embeddings defined by the tubular
neighborhood.  Since the moment map for the $\UU$-action is given in
this model by $\phi(x,t)=t$, the 1-form
$$\alpha:=-\iota\big(\d F(\f{\p}{\p t})\big)\om$$
is a principal connection for the $\UU$-action. Moreover $\alpha$ is 
adapted to $Y_\sig$ over $F(Y_\sig'\times (s-\eps,s+\eps))$.
 
The integral
$$ \rho_{t} := \int_s^t  F_{u}^* \alpha\ \d u\in \Om^1(\phi^{-1}(s))$$
is well-defined for $t\in (s-\eps,s+\eps)$, and is adapted to 
$Y_\sig$ over $F(Y_\sig'\times (s-\eps,s+\eps))$. 
The form $\rho_t$ satisfies
$$ F_{t}^* \om - F_s^* \om = \d \rho_{t} ,
\ \ \  F_{t}^* \Phi - F_s^*\Phi = -\rho_{t}^\sharp.$$
Since $\alpha$ is a principal connection on the tubular neighborhood,
$\rho_{t}/(t-s)$ is a principal connection for the $\UU$-action on
$\phinv(s)$, and the equivariant 2-form
$$\f{d\rho_{t}-2\pi i\,\rho_t^\sharp}{t-s}$$
represents its equivariant curvature. In particular it is  
basic, i.e descends to a closed equivariant 2-form $\gamma_t+2\pi
i\Psi_t$ on $X_s$, with 
$$ f_t^*{\om}_t=\om_s+(t-s)\gamma_s,$$
$$ f_t^*{\Phi}_t={\Phi}_s+(t-s)\Psi_s.$$
Since the cohomology class of $\gamma_t$ does not depend on $t$, we
have $\gamma_s=\gamma_t+\d\beta_t,\, \Psi_s=\Psi_t-\beta_t^\sharp$
where $\beta_t$ is given by the transgression formula.
\end{proof}

We now apply this result to the setting of Theorem 
\ref{ClassicalInductionProperty}. 

\begin{theorem}  \labell{DHApp}
Let $(M_1,\om_1)$ and $(M_2,\om_2)$ be compact
Hamiltonian $G$-manifolds. For $a>0$ let 
$M_i^{(a)}$ denote the same $G$-manifolds with re-scaled
2-forms $\om_i^{(a)} = a \om_i$.  Then 
there exists a family 
of equivariant diffeomorphisms, continuous and piecewise smooth in $a$,
$$ \phi_a:\, \Ind(M_1) \fus \Ind(M_2) 
\to  \Ind(M_1^{(a)}) \fus \Ind(M_2^{(a)}) $$
such that
$\phi_a^* \left( \Ind(\om_1^{(a)}) \fus
\Ind(\om_2^{(a)}) + 2\pi i \Ind( \Phi_1^{(a)}) \fus
\Ind( \Phi_2^{(a)}) \right)$ varies linearly in cohomology (in the sense 
of Proposition \ref{DH}).
\end{theorem}

\begin{proof}  For every $s=a^{-1}>0$ we define
$$ \Ind_s(M_i) := \Ind(M_i^{(1/s)})^{(s)} $$
which is a Hamiltonian $LG$-manifold at level $s$.  Equivalently, 
$\Ind_s(M_i)\cong LG\times_G M_i$ is the unique Hamiltonian 
$LG$-manifold with moment map at level $s$ for which the 2-form and moment map
restrict to the given ones on $M_i$. The 2-form on $\Ind_s(M_i)$ 
differs from the 2-form on $\Ind(M_i)$ by $(s-1)$ times the pull-back
of the symplectic form on $LG/G$.
Note that $\Ind(M_1^{(a)}) \fus \Ind(M_2^{(a)})$ is equivariantly
diffeomorphic to
$$ \Ind_s(M_1) \fus \Ind_s(M_2) := \M(\Sig^3_0)^{(s)} \times
\Ind_s(M_1) \times \Ind_s(M_2) \qu \diag(LG^2) $$
with symplectic forms equal after scalar multiplication by $s$.  We
show that there exists a Hamiltonian $\widehat{LG}$-manifold $W$ such
that the quotients $W_s$ with respect to the central $\UU$
are $\Ind_s(M_1) \fus \Ind_s(M_2)$, which
proves the Theorem by Proposition \ref{DH}.

To construct $W$ let $L \ra \M(\Sig^3_0) \times \Ind(M_1) \times
\Ind(M_2)$ be the pullback of the pre-quantum bundle on the base,
$$
L(\Sig^3_0) \boxtimes L(\Sig^1_0) \boxtimes L(\Sig^1_0).
$$
equipped with the pre-quantum connection.
Let $\alpha$ be the corresponding principal connection on the unit
circle bundle $U(L)$ and
$$ X = U(L) \times \R, \ \ \om_X = \pi^*\om + d(\phi-1,\alpha) $$
where $\om$ is the
closed two-form on $ \M(\Sig^3_0) \times \Ind(M_1) \times
\Ind(M_2)$ and $\phi$ is the coordinate on $\R$.
The quotient $X_s =
\phinv(s)/\UU$ is given by
$$ X_s \cong \M^s(\Sig^3_0) \times \Ind_s(M_1) \times
\Ind_s(M_2) .$$
Therefore if we set $ W = X \qu \diag(LG^2) $ then the central $\UU$
acts freely, and
$$ W_s = \Ind_s(M_1) \fus \Ind_s(M_2) $$
as claimed.  
\end{proof}

\section{Line bundles}

We now discuss extensions of the above results to include line
bundles, which will be required in the sequel \cite{me:co2} to 
this paper.

Let $M_1$ and $M_2$ be Hamiltonian $LG$-manifolds with proper moment
maps.  Given $\LG$-equivariant line bundles $L_1 \ra M_1, \ L_2 \ra
M_2$ at some level $k \in \Z$, we define%

$$L_1\fus L_2:=\big(L_1\boxtimes L_2\boxtimes
L(\Sig_0^3)^{\otimes k}\big)\qu\diag(LG^2) \ra M_1 \fus M_2. $$
If $L_1$ and $L_2$ are pre-quantum line bundles, then
$L_1\fus L_2$ is a pre-quantum bundle for $M_1\fus M_2$.

If $M$ is a Hamiltonian $G$-manifold with $G$-equivariant 
line bundle $L$, we define an $\widehat{LG}$-equivariant 
line bundle at level 
$\levi$ by 
$$ \Ind_\levi(L) := \widehat{LG} \times_{\widehat{G}} L \to \Ind(M)$$
where the central $U(1)$ acts on $L$ with fiber-weight $\levi$.  If
$L$ is a pre-quantum line bundle for $M$ then $\Ind_\levi(L)$ is a
pre-quantum line bundle for $\Ind_\levi(M)=\Ind(M^{(1/\levi)})^{(\levi)}$.
Note that the central extension $\widehat{G}=\widehat{LG}_{\{0\}}$
splits $\widehat{G}=G\times\UU$ because the restriction to $\g$ of the
defining cocycle is zero.

\begin{theorem}\label{Ext}
Let $(M_i,\omega_i,\Phi_i)$ $(i=1,2)$ be compact Hamiltonian
$G$-manifolds with $G$-equivariant line bundles $L_i \ra M_i$. Then
there exists an $LG$-equivariant diffeomorphism
$$
\phi:\, \Ind(M_1)\fus \Ind(M_2) \to \Ind(M_1\times M_2)
$$
such that the equivariant 2-forms on both sides are 
cohomologous (as in Theorem \ref{ClassicalInductionProperty}) 
and such that $\phi$ lifts to 
an equivariant isomorphism of line bundles:
$$ \hat{\phi}:\,\Ind_\levi(L_1) \fus \Ind_\levi(L_2) \to 
\Ind_\levi(L_1 \boxtimes L_2) .$$
\end{theorem}

\begin{proof} 
Let $\Psi:\,\M(\Sig_0^3)\to (L\g^*)^3$ be the moment map for
$\M(\Sig_0^3)$.  Recall from the proof of Theorem
\ref{ClassicalInductionProperty} that there exists a $G^2$-invariant
neighborhood $V\subset (U_{\{0\}})^2\subset (L\g^*)^2$ of $0$ such
that $\Psi^{-1}(V \times U_{\{0\}})$ is $G^3$-equivariantly
symplectomorphic to a neighborhood $U$ of the zero section in
$T^*G^2$.  First consider the special case when $\Phi_1(M_1) \times
\Phi_2(M_2) \subset V$.  Since equivariant line bundles over cotangent
bundles are trivial, the pullback of $L(\Sig^3_0)$ to $U$ is
equivariantly isomorphic to the trivial bundle $U\times\C$. Therefore
\begin{eqnarray*} 
 \Ind_\levi(L_1) \fus \Ind_\levi(L_2) &=& \Ind_\levi(L_1) \boxtimes
\Ind_\levi(L_2) \boxtimes L(\Sig^3_0)^{\otimes k} \qu LG^2 \\ &=&
\Ind_\levi(L_1) \boxtimes \Ind_\levi(L_2) \boxtimes \Ind_\levi(U \times \C )
\qu LG^2 \\ &=& \Ind_\levi(L_1 \boxtimes L_2 \boxtimes (U \times \C) \qu
G^2 ) \\ &=& \Ind_\levi(L_1 \boxtimes L_2).
\end{eqnarray*}

To reduce the general case to this special case, we rescale 
the 2-forms on $M_1,\,M_2$ by a factor $a>0$. Taking 
fusion products we obtain a family of $LG$-equivariant line bundles
$$ R^{(a)}=\Ind_\levi(L_1)\fus\Ind_\levi(L_2) \to \Ind(M_1^{(a)})\fus 
\Ind(M_2^{(a)}).$$
For $a$ small enough we are in the above situation and thus have an
isomorphism $R^{(a)}\cong \Ind_k(L_1\boxtimes L_2)$.  It remains to
show that the (continuous, piecewise smooth) family of diffeomorphisms
$\varphi_a$ given by Theorem \ref{DHApp} lifts to an equivariant
identification of line bundles $R^{(a)}$.  By continuity the
equivariant isomorphism class of the line bundle $\varphi_a^*R^{(a)}$
does not change, i.e. there are $LG$-equivariant 
isomorphisms $\varphi_a^*R^{(a)}\cong
\Ind_\levi(L_1)\fus\Ind_\levi(L_2)$ covering the identity map on
$\Ind(M_1)\fus \Ind(M_2)$.  This follows from the existence of an
$LG$-invariant connection on the family.  Such a connection may be
constructed using cross-sections and a partition of unity as in Lemma
\ref{Adapted}.  Alternatively, the claim follows because the family of
line bundles $(\varphi_a^*R^{(a)})\otimes
\big(\varphi_1^*R^{(1)} \big)^{-1}$ descends to a family of 
$G$-equivariant line bundles on the holonomy manifold 
$\Hol\big(\Ind(M_1)\fus \Ind(M_2)\big)\cong M_1\times M_2$, equal 
to the trivial bundle for $a=1$.
\end{proof}






\end{document}